\documentclass[english,12pt]{article}

\pdfoutput=1

\usepackage{amssymb}
\usepackage{amsthm}
\usepackage{amsmath}
\usepackage{amsrefs}
\usepackage{array}
\usepackage{enumerate}
\usepackage{mathdots}
\usepackage{multicol}
\usepackage{algorithm}
\usepackage[noend]{algpseudocode}
\usepackage{calc}
\usepackage{tikz}
\usepackage{tikz-qtree}
\usetikzlibrary{calc}
\usetikzlibrary{shapes}

\newtheorem{defin}{Definition}
\newtheorem{prop}{Proposition}
\newtheorem{theorem}{Theorem}
\newtheorem{lemma}{Lemma}
\newtheorem{cor}{Corollary}

\newcommand{\pb}[1]{\left(#1\right)}
\newcommand{\p}[1]{(#1)}
\newcommand{\st}[1]{\left\{#1\right\}}

\newcommand{\bk}[1]{\left[#1\right]}

\newcommand{\fl}[1]{\left\lfloor#1\right\rfloor}

\newcommand{\quot}[1]{``#1''}

\newcommand{\fwd}{\p{\Longrightarrow}}
\newcommand{\bwd}{\p{\Longleftarrow}}

\newcommand{\iffpf}[2]{
\begin{description}
\item[$\fwd$] #1
\item[$\bwd$] #2
\end{description}
}

\newcommand{\limf}[3]{\lim_{#1\rightarrow#2}{#3}}

\newcommand{\limi}[2]{\limf{#1}{\infty}{#2}}

\renewcommand{\bk}[1]{[#1]}

\newcommand{\seq}{\pb}

\newcommand{\NN}{L}

\newcommand{\ninduct}{m}

\newcommand{\dseqend}[4]{#1_{#2},#1_{#2-1},\ldots,#1_{#3},#1_{#4}}
\newcommand{\dseq}[2]{\dseqend{#1}{#2}{1}{0}}

\newcommand{\typel}[1]{Type~$#1$}

\newcounter{temthm}
\newcounter{temthma}
\newcommand{\mainthm}{Let $a=\bk{\lambda_1,\lambda_2,\ldots,\lambda_k}$ be a linear-recurrent sequence, and let $s$ be a nonnegative integer.  Let $T=T\bk{\lambda_1,\lambda_2,\ldots,\lambda_k}$, and let $t$ denote the minimum index such that $L_{T,s}\p{t}=a_k$. Define the nested recurrence
\[
C_a\p{n}=\sum_{\ell=0}^{\Lambda-1}C_a\p{R\p{C_a, n,\mu_\ell, s, \ell}}.
\]
The recurrence $C_a$ generates the sequence $\seq{L_{T,s}\p{n}}_{n\geq1}$ using terms $L_{T,s}\p{1}$ through $L_{T,s}\p{t}$ as initial conditions.}% TODO negative $s$ changes length of IC

\begin{document}
%
%\lstset{breaklines=true, morecomment=[l]{//}, frame=single, showstringspaces=false, numbers=left}
%
%\begin{flushright}
%Nathan Fox\\
%Math 640\\
%\end{flushright}
%\begin{center}
%\LARGE{On Aperiodic Subtraction Games with Bounded Nim Sequence}
%\end{center}
%\namehead{Math 640}{\LARGE{On Aperiodic Subtraction Games with Bounded Nim Sequence}}
\title{Connecting Slow Solutions to Nested Recurrences with Linear Recurrent Sequences}
\author{Nathan Fox\footnote{Mathematical and Statistics, Canisius College, Buffalo, New York,
\texttt{fox42@canisius.edu}
}}
\date{}

\maketitle

\begin{abstract}
Labeled infinite trees provide combinatorial interpretations for many integer sequences generated by nested recurrence relations. Typically, such sequences are monotone increasing. Several of these sequences also have straightforward descriptions in terms of how often each value in the sequence occurs. In this paper, we generalize the most classical examples to a larger family of sequences parametrized by linear recurrence relations. Each of our sequences can be constructed in three different ways: via a nested recurrence relation, from labeled infinite trees, or by using Zeckendorf-like strings of digits to describe its frequency sequence. We conclude the paper by discussing the asymptotic behaviors of our sequences.
\end{abstract}

\section{Introduction}
Nested recurrence relations provide a concise method of generating a wide variety of sequences of integers. Perhaps the earliest and most well-known nested recurrence relation is Hofstadter's $Q$-recurrence~\cite{geb}
\[
Q\p{n}=Q\p{n-Q\p{n-1}}+Q\p{n-Q\p{n-2}}.
\]
Hofstadter considered this recurrence with initial conditions $Q\p{1}=Q\p{2}=1$, calling the result a \emph{meta-Fibonacci sequence}. This sequence, the Hofstadter $Q$-sequence, behaves chaotically, though with a hint of structure. Most questions about this sequence remain open, including whether it is well-defined for all positive indices. Computational evidence~\cites{pinn, oeis} suggests that it is and that $Q\p{n}$ remains close to $\frac{n}{2}$.

While the $Q$-sequence itself has remained out of grasp, several results have been proved about related sequences. One approach focuses on choosing a nested recurrence and finding initial conditions that result in predictable behavior~\cites{golomb, rusk, genrusk, gengol, symbhof,hof1thruN}. Often, the resulting structure is eventually an interleaving of well understood sequences. Another approach focuses on finding nested recurrence relations where simple initial conditions, often several $1$'s, appear to result in a monotone increasing sequence~\cites{golomb,isgurthesis,slowtrihof,erickson,hofv,tanny,con,mallows}. Our results come about as an extension of this second approach.

Of particular interest among monotone sequences arising from nested recurrences are those in which every positive integer appears, so-called \emph{slow sequences}~\cite{isgur2}. Equivalently, a sequence of integers is slow if its first term is $1$ and successive differences between successive terms are always $0$ or $1$. A slow sequence can be described by its \emph{frequency sequence}, a sequence of positive integers $f_1,f_2,\ldots$, where $f_i$ enumerates the number of occurrences of the integer $i$.

Many slow sequences arising from nested recurrence relations have combinatorial interpretations related to counting leaves in infinite labeled trees~\cites{isgur1,isgur2,erickson,isgurthesis,ruskeynegative}. These trees are constructed recursively; we explain the construction and labeling process in Section~\ref{ss:tree}. Classically, sequences arising from trees are necessarily slow, though some recent variations on the tree-based methodology have allowed generating non-slow sequences from trees as well~\cites{nonslow,nonmonotgol}. Such combinatorial interpretations often lead directly to an asymptotic description of the sequence.

Three of the most well-known slow sequences with tree-based combinatorial interpretations are the Conolly sequence~\cite{con}, the Tanny sequence~\cite{tanny}, and the Golomb sequence~\cite{golomb}. Tables~\ref{tab:con},~\ref{tab:tan}, and~\ref{tab:gol} provides some basic information about each of these sequences. The Golomb recurrence is beyond the scope of our work, so we henceforth set it aside. The structures of the Conolly and Tanny sequences, on the other hand, are quite similar to one another. The frequency of the number $N$ in the Conolly sequence is the smallest positive integer $t$ such that $2^t$ does not divide $N$~\cite{oeis}. For the Tanny sequence, the frequency sequence is nearly identical, except when $N$ is power of $2$. In those cases, the Tanny sequence contains one more $N$ term than the Conolly sequence~\cite{oeis}. These frequency counts and the relationship between these two sequences follows from their related combinatorial interpretations~\cite{isgur2}.%plus the number of zeroes at the end of the binary representation of $N$~\cite{oeis}.

\begin{table}
\begin{center}
\begin{tabular}{|l||l|}\hline
\textbf{Recurrence} & $C\p{n}=C\p{n-C\p{n-1}}+C\p{n-1-C\p{n-2}}$ \\\hline%& $T\p{n}=T\p{n-1-T\p{n-1}}+T\p{n-2-T\p{n-2}}$ & $G\p{n}=G\p{n-G\p{n-1}}+1$\\\hline
\textbf{Initial Conditions}\footnotemark &  $C\p{1}=1$, $C\p{2}=2$\\\hline% & $T\p{1}=T\p{2}=T\p{3}=1$ & $G\p{1}=1$\\\hline
\textbf{OEIS} & A046699 \\\hline%& A006949 & A002024\\\hline
\textbf{Terms}  & $1,2,2,3,4,4,4,5,6,6,7,8,8,8,8,9,\ldots$ \\\hline%& $1,1,1,2,2,2,3,4,4,4,4,5,6,6,7,\ldots$ & $1,2,2,3,3,3,4,4,4,4,5,5,5,5,5,\ldots$\\\hline
\textbf{Frequency Sequence} & $1,2,1,3,1,2,1,4,1,2,1,3,1,2,1,5,1,2,1,3,\ldots$\\\hline
\textbf{Asymptotic} & $C\p{n}\sim\frac{n}{2}$\\\hline
\end{tabular}
\end{center}
\caption{Information about the Conolly sequence}
\label{tab:con}
\end{table}
\footnotetext{\label{fn:ic}The classical initial conditions are $C\p{1}=C\p{2}=1$ for the Conolly sequence and $T\p{1}=T\p{2}=T\p{3}=1$ for the Tanny sequence. Our initial conditions result in the same sequences, both with an initial~$1$ removed. This leads to cleaner combinatorial interpretations for both.}

\begin{table}
\begin{center}
\begin{tabular}{|l||l|}\hline
\textbf{Recurrence} & $T\p{n}=T\p{n-1-T\p{n-1}}+T\p{n-2-T\p{n-2}}$\\\hline% & $G\p{n}=G\p{n-G\p{n-1}}+1$\\\hline
\textbf{Initial Conditions}\textsuperscript{\ref{fn:ic}} &  $T\p{1}=T\p{2}=1$, $T\p{3}=2$\\\hline% & $G\p{1}=1$\\\hline
\textbf{OEIS} & A006949 \\\hline%& A002024\\\hline
\textbf{Terms}  & $1,1,2,2,2,3,4,4,4,4,5,6,6,7,8,8,8,8,8,9\ldots$ \\\hline%& $1,2,2,3,3,3,4,4,4,4,5,5,5,5,5,\ldots$\\\hline
\textbf{Frequency Sequence} & $2,3,1,4,1,2,1,5,1,2,1,3,1,2,1,6,1,2,1,3,\ldots$\\\hline
\textbf{Asymptotic} & $T\p{n}\sim\frac{n}{2}$\\\hline
\end{tabular}
\end{center}
\caption{Information about the Tanny sequence}
\label{tab:tan}
\end{table}

\begin{table}
\begin{center}
\begin{tabular}{|l||l|}\hline
\textbf{Recurrence} & $G\p{n}=G\p{n-G\p{n-1}}+1$\\\hline
\textbf{Initial Conditions} & $G\p{1}=1$\\\hline
\textbf{OEIS} & A002024\\\hline
\textbf{Terms}  & $1,2,2,3,3,3,4,4,4,4,5,5,5,5,5,6,6,6,6,6,6,\ldots$\\\hline
\textbf{Frequency Sequence} & $1,2,3,4,5,6,7,8,9,10,11,12,13,14,15,16,\ldots$\\\hline
\textbf{Asymptotic} & $G\p{n}\sim\sqrt{2n}$\\\hline
\end{tabular}
\end{center}
\caption{Information about the Golomb sequence}
\label{tab:gol}
\end{table}

In this paper, we generalize the Conolly and Tanny sequences in a new way. These sequences have an intimate connection to the powers of $2$, and generalizations are known for powers of any larger integer~\cites{isgur2,ruskeynegative}. We generalize sequences of powers to sequences generated by linear recurrence relations with positive integer coefficients, such as the Fibonacci sequence. Such sequences can be used to generate frequency sequences of slow sequences, and we show how these slow sequences can be obtained from nested recurrence relations, possibly with several layers of nesting (Section~\ref{sec:treerec}). Then, we provide an alternative combinatorial interpretation in terms of strings, a generalization of Zeckendorf representations of integers (Section~\ref{sec:genzeck}), and we use that to derive asymptotics for all of our sequences (Section~\ref{sec:asympt}). %Once we have all of the pieces, we discuss three examples in detail (Section~\ref{sec:examples}). 
Finally, we conclude by describing some open problems and future directions (Section~\ref{sec:fut}).

\subsection{An Overview of the Tree Methodology}\label{ss:tree}
To illustrate the tree methodology, we show how to apply it to the family of recurrences $C_s\p{n}=C_s\p{n-s-C_s\p{n-1}}+C_s\p{n-s-1-C_s\p{n-2}}$, where $s$ is an integer\footnote{Several authors prefer the convention where $C_s\p{n}=C_s\p{n-s+1-C_s\p{n-1}}+C_s\p{n-s-C_s\p{n-2}}$, adding $1$ to each of our $s$ values.}  The situation when $s<0$ is more complicated; we henceforth assume $s\geq0$ unless explicitly specified otherwise. When $s\geq0$, we consider an initial condition of $s+1$ $1$'s followed by a $2$. The case $s=0$ is the Conolly sequence~\cite{con}, and the case $s=1$ is the Tanny sequence~\cite{tanny}.  Henceforth, we call members of this family \emph{generalized Conolly sequences}.  A single tree structure suffices to analyze all of these sequences simultaneously. What follows parallels a discussion found in~\cite{isgur2} and~\cite{ruskeynegative}.

We begin by defining a \emph{skeleton}: an infinite tree that we later add labels to.  The precise form of the skeleton is determined by the recurrence, %~\cite{isgurthesis},
but it is always defined starting from its leftmost leaf and working upward and rightward.  Skeletons are built of \emph{levels}; the leaves are at level~$0$, and the level of a non-leaf is its height above its leaf children.  In particular, skeletons are \emph{perfect}; all leaves are at the same level, so the level of every node is well-defined.  The nodes on each level have an ordering from left to right.  We will find the following definitions useful going forward:
\begin{defin}\label{def:type}
The \emph{type} of a node in a skeleton is the number of left siblings it has.
\end{defin}
\begin{defin}\label{def:special}
A node is \emph{special} if it and all of its ancestors are \typel{0}. 
\end{defin}
Equivalently, a node is special if it is the leftmost node on its level; each level has exactly one special node.  The term \emph{$i^{th}$ special node} refers to the special node on level~$i$.
\begin{defin}\label{def:subtree}
A \emph{subtree} of a skeleton is a node along with all of its descendants. If the root of a subtree is special, we call the subtree a \emph{left subtree}.
\end{defin}
\begin{defin}\label{def:stree}
The \emph{$i^{th}$ sub-skeleton} of a skeleton is the tree rooted at the $i^{th}$ special node and containing all of its non-special children and all of their descendants.
\end{defin}
Note that the nodes of a skeleton are the disjoint union of the nodes of its sub-skeletons, and sub-skeletons $0$ through $i$ combine to give the 
left subtree rooted at the $i^{th}$ special node.
%tree rooted at the $i^{th}$ special node containing all of its descendants.

In addition to each node having a type, nodes in a skeleton are classified into one of three categories:
\begin{description}
\item[Leaves:] Any leaf node.  We draw leaves as ellipses.
\item[Supernodes:] Any special node that is not a leaf.  (There is only one special leaf: the leftmost leaf.)  The term \emph{$i^{th}$ supernode} will be used to refer to the supernode on level~$i$.  We draw supernodes as rectangles.
\item[Regular Nodes:] Any node that is neither a leaf nor a supernode.  We draw regular nodes as circles.
\end{description}

In the case of the generalized Conolly sequences, the skeleton is a full binary tree (Figure~\ref{fig:conskel}).  
\begin{figure}
\begin{center}
\resizebox{\textwidth-0.5in}{!}{
\tikzstyle{super}=[rectangle, draw]
\tikzstyle{regular}=[circle, draw]
\tikzstyle{leaf}=[ellipse, draw]
\begin{tikzpicture}[level distance = 2.2cm, grow=down]
\Tree[.{$\phantom{aaa}\iddots$}
[.\node[super]{$\phantom{0000}$};
  [.\node[super]{$\phantom{0000}$};
    [.\node[super]{$\phantom{0000}$};
      [.\node[super]{$\phantom{0000}$};
        [.\node[leaf]{$\phantom{0000}$};
        ]
        [.\node[leaf]{$\phantom{0000}$};
        ]
      ]
      [.\node[regular]{$\phantom{0000}$};
        [.\node[leaf]{$\phantom{0000}$};
        ]
        [.\node[leaf]{$\phantom{0000}$};
        ]
      ]
    ]
    [.\node[regular]{$\phantom{0000}$};
      [.\node[regular]{$\phantom{0000}$};
        [.\node[leaf]{$\phantom{0000}$};
        ]
        [.\node[leaf]{$\phantom{0000}$};
        ]
      ]
      [.\node[regular]{$\phantom{0000}$};
        [.\node[leaf]{$\phantom{0000}$};
        ]
        [.\node[leaf]{$\phantom{0000}$};
        ]
      ]
    ]
  ]
  [.\node[regular]{$\phantom{0000}$};
    [.\node[regular]{$\phantom{0000}$};
      [.\node[regular]{$\phantom{0000}$};
        [.\node[leaf]{$\phantom{0000}$};
        ]
        [.\node[leaf]{$\phantom{0000}$};
        ]
      ]
      [.\node[regular]{$\phantom{0000}$};
        [.\node[leaf]{$\phantom{0000}$};
        ]
        [.\node[leaf]{$\phantom{0000}$};
        ]
      ]
    ]
    [.\node[regular]{$\phantom{0000}$};
      [.\node[regular]{$\phantom{0000}$};
        [.\node[leaf]{$\phantom{0000}$};
        ]
        [.\node[leaf]{$\phantom{0000}$};
        ]
      ]
      [.\node[regular]{$\phantom{0000}$};
        [.\node[leaf]{$\phantom{0000}$};
        ]
        [.\node[leaf]{$\phantom{0000}$};
        ]
      ]
    ]
  ]
]
\edge[draw=none]; {} ]
\end{tikzpicture}
}
\end{center}
\caption{The lower left portion of the skeleton for the generalized Conolly sequences, up to the fourth supernode}
\label{fig:conskel}
\end{figure}
A \emph{preorder traversal} of a skeleton is similar to a preorder traversal of a rooted tree, except that a skeleton has no root.  Instead, we traverse the sub-skeletons (Definition~\ref{def:stree}) in increasing order, and we traverse each sub-skeleton in preorder.  We use preorder traversal to generate a \emph{labeled skeleton}.  For the recurrence $C_s$, as we visit each node, we write the smallest unused positive integer in it if it is a leaf or a regular node, and we write the smallest $s$ unused positive integers in it if it is a supernode (Figure~\ref{fig:conskelab}).  For a general skeleton $T$ and an integer parameter $s\geq0$, we obtain the labeled skeleton, denoted $T^{(s)}$, via a preorder traversal, putting $s$ labels in each supernode and one label in each leaf or regular node. In the cases where $s<0$, zero labels are placed in each supernode, and the leftmost $-s$ regular nodes on each level are left empty as well~\cites{isgur2,ruskeynegative}.
\begin{figure}
\begin{center}
\resizebox{\textwidth-0.5in}{!}{
\tikzstyle{super}=[rectangle, draw]
\tikzstyle{regular}=[circle, draw]
\tikzstyle{leaf}=[ellipse, draw]
\begin{tikzpicture}[level distance = 1.5cm, grow=down]
\Tree[.{$\phantom{aaa}\iddots$}
[.\node[super]{$19,20$};
  [.\node[super]{$10,11$};
    [.\node[super]{$\phantom{0}5,6\phantom{0}$};
      [.\node[super]{$\phantom{0}2,3\phantom{0}$};
        [.\node[leaf]{$1\phantom{0}$};
        ]
        [.\node[leaf]{$4\phantom{0}$};
        ]
      ]
      [.\node[regular]{$7\phantom{0}$};
        [.\node[leaf]{$8\phantom{0}$};
        ]
        [.\node[leaf]{$9\phantom{0}$};
        ]
      ]
    ]
    [.\node[regular]{$12$};
      [.\node[regular]{$13$};
        [.\node[leaf]{$14$};
        ]
        [.\node[leaf]{$15$};
        ]
      ]
      [.\node[regular]{$16$};
        [.\node[leaf]{$17$};
        ]
        [.\node[leaf]{$18$};
        ]
      ]
    ]
  ]
  [.\node[regular]{$21$};
    [.\node[regular]{$22$};
      [.\node[regular]{$23$};
        [.\node[leaf]{$24$};
        ]
        [.\node[leaf]{$25$};
        ]
      ]
      [.\node[regular]{$26$};
        [.\node[leaf]{$27$};
        ]
        [.\node[leaf]{$28$};
        ]
      ]
    ]
    [.\node[regular]{$29$};
      [.\node[regular]{$30$};
        [.\node[leaf]{$31$};
        ]
        [.\node[leaf]{$32$};
        ]
      ]
      [.\node[regular]{$33$};
        [.\node[leaf]{$34$};
        ]
        [.\node[leaf]{$35$};
        ]
      ]
    ]
  ]
]
\edge[draw=none]; {} ]
\end{tikzpicture}
}
\end{center}
\caption{The lower left portion of the labeled skeleton for the sequence $C_2\p{n}$, up to the fourth supernode}
\label{fig:conskelab}
\end{figure}

%Leaf-counting function
From a skeleton $T$ and its corresponding labeled skeleton $T^{(s)}$ for some value of $s$, we can obtain its \emph{leaf-counting function} $L_{T,s}\p{n}$.  For positive integer $n$, we define $L_{T,s}\p{n}$ to be the number of leaves of $T^{(s)}$ with labels at most $n$.  For example, if $T$ is the skeleton in Figure~\ref{fig:conskel}, Figure~\ref{fig:conskelab} shows $T^{(2)}$ and that $L_{T,2}\p{8}=3$ and $L_{T,2}\p{30}=12$.  Note that for any labeled skeleton, the sequence $\seq{L_{T,s}\p{n}}_{n\geq1}$ is slow.

\begin{figure}
\begin{center}
\resizebox{\textwidth-0.5in}{!}{
%drawTree(s = 2, n = 19)
\tikzstyle{super}=[rectangle, draw]
\tikzstyle{regular}=[circle, draw]
\tikzstyle{leaf}=[ellipse, draw]
\begin{tikzpicture}[level distance = 1.5cm, grow=down]
\Tree[.{$\phantom{aaa}\iddots$}
[.\node[super]{$\phantom{0}19\phantom{00}$};
  [.\node[super]{$10,11$};
    [.\node[super]{$\phantom{0}5,6\phantom{0}$};
      [.\node[super]{$\phantom{0}2,3\phantom{0}$};
        [.\node[leaf]{$1\phantom{0}$};
        ]
        [.\node[leaf]{$4\phantom{0}$};
        ]
      ]
      [.\node[regular]{$7\phantom{0}$};
        [.\node[leaf]{$8\phantom{0}$};
        ]
        [.\node[leaf]{$9\phantom{0}$};
        ]
      ]
    ]
    [.\node[regular]{$12$};
      [.\node[regular]{$13$};
        [.\node[leaf]{$14$};
        ]
        [.\node[leaf]{$15$};
        ]
      ]
      [.\node[regular]{$16$};
        [.\node[leaf]{$17$};
        ]
        [.\node[leaf]{$18$};
        ]
      ]
    ]
  ]
  [.\node[regular]{$\phantom{00}$};
    [.\node[regular]{$\phantom{00}$};
      [.\node[regular]{$\phantom{00}$};
        [.\node[leaf]{$\phantom{00}$};
        ]
        [.\node[leaf]{$\phantom{00}$};
        ]
      ]
      [.\node[regular]{$\phantom{00}$};
        [.\node[leaf]{$\phantom{00}$};
        ]
        [.\node[leaf]{$\phantom{00}$};
        ]
      ]
    ]
    [.\node[regular]{$\phantom{00}$};
      [.\node[regular]{$\phantom{00}$};
        [.\node[leaf]{$\phantom{00}$};
        ]
        [.\node[leaf]{$\phantom{00}$};
        ]
      ]
      [.\node[regular]{$\phantom{00}$};
        [.\node[leaf]{$\phantom{00}$};
        ]
        [.\node[leaf]{$\phantom{00}$};
        ]
      ]
    ]
  ]
]
\edge[draw=none]; {} ]
\end{tikzpicture}
}
\end{center}
\caption{$T^{(2)}\p{19}$ of the skeleton for the generalized Conolly sequences, up to the fourth supernode}
\label{fig:conskeln}
\end{figure}
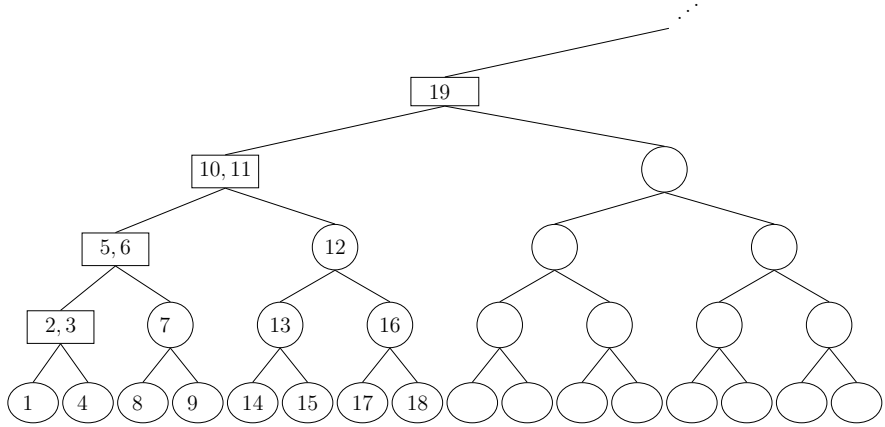

%T, T_s(n)
As an additional useful piece of notation, let $T^{(s)}\p{n}$ denote the skeleton $T$ with labels $1$ through $n$ placed in the tree according to the rule for sequence $C_s$ (Figure~\ref{fig:conskeln}).  
The key property is the following:
\begin{theorem}\cite{ruskeynegative}\label{thm:contree}
Let $T$ be the full binary tree skeleton (Figure~\ref{fig:conskel}).  For all $n\geq1$ and $s\geq0$, $L_{T,s}\p{n}=C_s\p{n}$.
\end{theorem}
The proof illustrates how the structure of labeled skeletons is related to nested recurrences.
%Pruning
\begin{proof}
We argue by induction on $n$.  First, note that the lower left leaf has the label $1$, the first supernode has the labels $2$ through $2+s-1$ (or no labels if $s=0$), and the second leaf has label $s+2$.  This means that the first $s+1$ terms of the sequence $L_{T,s}\p{n}$ are $1$, and term $s+2$ is $2$.  These terms match the initial conditions of $C_s\p{n}$.

Now, suppose $n>s+2$, and suppose that $C_s\p{\ninduct}=L_{T,s}\p{\ninduct}$ for all $\ninduct<n$.  
Note that the value of $L_{T,s}\p{n}$ equals the number of labeled leaves in $T^{(s)}\p{n}$.  Now, define $L_{T,s,0}\p{n}$ to be the number of labeled \typel{0} leaves in $T^{(s)}\p{n}$, and define $L_{T,s,1}\p{n}$ analogously for \typel{1} leaves.  Since the tree $T$ is binary, every node is either \typel{0} or \typel{1}, so
\[
L_{T,s}\p{n}=L_{T,s,0}\p{n}+L_{T,s,1}\p{n}.
\]
To show that $L_{T,s}\p{n}=C_s\p{n}$, we will show that $L_{T,s,0}\p{n}=L_{T,s}\p{n-s-L_{T,s}\p{n-1}}$ and that $L_{T,s,1}\p{n}=L_{T,s}\p{n-s-1-L_{T,s}\p{n-2}}$.  This will show that the sequence $\seq{L_{T,s}\p{n}}_{n\geq1}$ satisfies the same recurrence as the sequence $\seq{C_s\p{n}}_{n\geq1}$.

%We now claim that $L_{T,s,\ell}\p{n}=C_a\p{R\p{C_a, \mu_\ell, s, \ell}}$.
The proof of the preceding fact comes from combinatorial manipulation of the tree.  We define a \emph{pruning} process for $T^{(s)}\p{n}$ consisting of the following five steps:
\begin{description}
\item[Deletion Step:] Delete all leaf labels in $T^{(s)}\p{n}$ that are less than $n$.
\item[Correction Step:] Delete all $s$ labels from the first supernode.
\item[Lifting Step:] Move the label $n$ to the first supernode.
\item[Restructuring Step:] Delete all the leaf nodes, thereby decreasing the level of each remaining node by $1$.  As a consequence, what was the first supernode becomes a leaf node, and all the regular nodes that were on level~$1$ also become leaf nodes.
\item[Relabeling Step:] Relabel the current labels by a preorder labeling.
\end{description}
%(See Figures~\ref{fig:con27},~\ref{fig:condel}, and~\ref{fig:conrel}.) 
See Figure~\ref{fig:condel} for the example of pruning $T^{(2)}\p{27}$. 
%\begin{figure}
%\begin{center}
%\includegraphics[width=400pt]{contree27.eps}
%\end{center}
%\caption{$T_2\p{27}$, prior to pruning}
%\label{fig:con27}
%\end{figure}
\begin{figure}
\begin{center}
\begin{tabular}{c|c}
\resizebox{\textwidth*\real{0.5}-0.2in}{!}{
%drawTree(s = 2, n = 27, height = 2.000000)
\tikzstyle{super}=[rectangle, draw]
\tikzstyle{regular}=[circle, draw]
\tikzstyle{leaf}=[ellipse, draw]
\begin{tikzpicture}[level distance = 2.0cm, grow=down]
\Tree[.{$\phantom{aaa}\iddots$}
[.\node[super]{$19,20$};
  [.\node[super]{$10,11$};
    [.\node[super]{$\phantom{0}5,6\phantom{0}$};
      [.\node[super]{$\phantom{0}2,3\phantom{0}$};
        [.\node[leaf]{$1\phantom{0}$};
        ]
        [.\node[leaf]{$4\phantom{0}$};
        ]
      ]
      [.\node[regular]{$7\phantom{0}$};
        [.\node[leaf]{$8\phantom{0}$};
        ]
        [.\node[leaf]{$9\phantom{0}$};
        ]
      ]
    ]
    [.\node[regular]{$12$};
      [.\node[regular]{$13$};
        [.\node[leaf]{$14$};
        ]
        [.\node[leaf]{$15$};
        ]
      ]
      [.\node[regular]{$16$};
        [.\node[leaf]{$17$};
        ]
        [.\node[leaf]{$18$};
        ]
      ]
    ]
  ]
  [.\node[regular]{$21$};
    [.\node[regular]{$22$};
      [.\node[regular]{$23$};
        [.\node[leaf]{$24$};
        ]
        [.\node[leaf]{$25$};
        ]
      ]
      [.\node[regular]{$26$};
        [.\node[leaf]{$27$};
        ]
        [.\node[leaf]{$\phantom{00}$};
        ]
      ]
    ]
    [.\node[regular]{$\phantom{00}$};
      [.\node[regular]{$\phantom{00}$};
        [.\node[leaf]{$\phantom{00}$};
        ]
        [.\node[leaf]{$\phantom{00}$};
        ]
      ]
      [.\node[regular]{$\phantom{00}$};
        [.\node[leaf]{$\phantom{00}$};
        ]
        [.\node[leaf]{$\phantom{00}$};
        ]
      ]
    ]
  ]
]
\edge[draw=none]; {} ]
\end{tikzpicture}
} & 
\resizebox{\textwidth*\real{0.5}-0.2in}{!}{
%drawTree(s = 2, n = 27, height = 2.000000, labels = ['', 2, 3, '', 5, 6, 7, '', '', 10, 11, 12, 13, '', '', 16, '', '', 19, 20, 21, 22, 23, '', '', 26, 27])
\tikzstyle{super}=[rectangle, draw]
\tikzstyle{regular}=[circle, draw]
\tikzstyle{leaf}=[ellipse, draw]
\begin{tikzpicture}[level distance = 2.0cm, grow=down]
\Tree[.{$\phantom{aaa}\iddots$}
[.\node[super]{$19,20$};
  [.\node[super]{$10,11$};
    [.\node[super]{$\phantom{0}5,6\phantom{0}$};
      [.\node[super]{$\phantom{0}2,3\phantom{0}$};
        [.\node[leaf]{$\phantom{00}$};
        ]
        [.\node[leaf]{$\phantom{00}$};
        ]
      ]
      [.\node[regular]{$7\phantom{0}$};
        [.\node[leaf]{$\phantom{00}$};
        ]
        [.\node[leaf]{$\phantom{00}$};
        ]
      ]
    ]
    [.\node[regular]{$12$};
      [.\node[regular]{$13$};
        [.\node[leaf]{$\phantom{00}$};
        ]
        [.\node[leaf]{$\phantom{00}$};
        ]
      ]
      [.\node[regular]{$16$};
        [.\node[leaf]{$\phantom{00}$};
        ]
        [.\node[leaf]{$\phantom{00}$};
        ]
      ]
    ]
  ]
  [.\node[regular]{$21$};
    [.\node[regular]{$22$};
      [.\node[regular]{$23$};
        [.\node[leaf]{$\phantom{00}$};
        ]
        [.\node[leaf]{$\phantom{00}$};
        ]
      ]
      [.\node[regular]{$26$};
        [.\node[leaf]{$27$};
        ]
        [.\node[leaf]{$\phantom{00}$};
        ]
      ]
    ]
    [.\node[regular]{$\phantom{00}$};
      [.\node[regular]{$\phantom{00}$};
        [.\node[leaf]{$\phantom{00}$};
        ]
        [.\node[leaf]{$\phantom{00}$};
        ]
      ]
      [.\node[regular]{$\phantom{00}$};
        [.\node[leaf]{$\phantom{00}$};
        ]
        [.\node[leaf]{$\phantom{00}$};
        ]
      ]
    ]
  ]
]
\edge[draw=none]; {} ]
\end{tikzpicture}
}\\\hline
\resizebox{\textwidth*\real{0.5}-0.2in}{!}{
%drawTree(s = 2, n = 27, height = 2.000000, labels = ['', '', '', '', 5, 6, 7, '', '', 10, 11, 12, 13, '', '', 16, '', '', 19, 20, 21, 22, 23, '', '', 26, 27])
\tikzstyle{super}=[rectangle, draw]
\tikzstyle{regular}=[circle, draw]
\tikzstyle{leaf}=[ellipse, draw]
\begin{tikzpicture}[level distance = 2.0cm, grow=down]
\Tree[.{$\phantom{aaa}\iddots$}
[.\node[super]{$19,20$};
  [.\node[super]{$10,11$};
    [.\node[super]{$\phantom{0}5,6\phantom{0}$};
      [.\node[super]{$\phantom{00000}$};
        [.\node[leaf]{$\phantom{00}$};
        ]
        [.\node[leaf]{$\phantom{00}$};
        ]
      ]
      [.\node[regular]{$7\phantom{0}$};
        [.\node[leaf]{$\phantom{00}$};
        ]
        [.\node[leaf]{$\phantom{00}$};
        ]
      ]
    ]
    [.\node[regular]{$12$};
      [.\node[regular]{$13$};
        [.\node[leaf]{$\phantom{00}$};
        ]
        [.\node[leaf]{$\phantom{00}$};
        ]
      ]
      [.\node[regular]{$16$};
        [.\node[leaf]{$\phantom{00}$};
        ]
        [.\node[leaf]{$\phantom{00}$};
        ]
      ]
    ]
  ]
  [.\node[regular]{$21$};
    [.\node[regular]{$22$};
      [.\node[regular]{$23$};
        [.\node[leaf]{$\phantom{00}$};
        ]
        [.\node[leaf]{$\phantom{00}$};
        ]
      ]
      [.\node[regular]{$26$};
        [.\node[leaf]{$27$};
        ]
        [.\node[leaf]{$\phantom{00}$};
        ]
      ]
    ]
    [.\node[regular]{$\phantom{00}$};
      [.\node[regular]{$\phantom{00}$};
        [.\node[leaf]{$\phantom{00}$};
        ]
        [.\node[leaf]{$\phantom{00}$};
        ]
      ]
      [.\node[regular]{$\phantom{00}$};
        [.\node[leaf]{$\phantom{00}$};
        ]
        [.\node[leaf]{$\phantom{00}$};
        ]
      ]
    ]
  ]
]
\edge[draw=none]; {} ]
\end{tikzpicture}
} & 
\resizebox{\textwidth*\real{0.5}-0.2in}{!}{
%drawTree(s = 2, n = 27, height = 2.000000, labels = ['', 27, '', '', 5, 6, 7, '', '', 10, 11, 12, 13, '', '', 16, '', '', 19, 20, 21, 22, 23, '', '', 26, ''])
\tikzstyle{super}=[rectangle, draw]
\tikzstyle{regular}=[circle, draw]
\tikzstyle{leaf}=[ellipse, draw]
\begin{tikzpicture}[level distance = 2.0cm, grow=down]
\Tree[.{$\phantom{aaa}\iddots$}
[.\node[super]{$19,20$};
  [.\node[super]{$10,11$};
    [.\node[super]{$\phantom{0}5,6\phantom{0}$};
      [.\node[super]{$\phantom{0}27\phantom{00}$};
        [.\node[leaf]{$\phantom{00}$};
        ]
        [.\node[leaf]{$\phantom{00}$};
        ]
      ]
      [.\node[regular]{$7\phantom{0}$};
        [.\node[leaf]{$\phantom{00}$};
        ]
        [.\node[leaf]{$\phantom{00}$};
        ]
      ]
    ]
    [.\node[regular]{$12$};
      [.\node[regular]{$13$};
        [.\node[leaf]{$\phantom{00}$};
        ]
        [.\node[leaf]{$\phantom{00}$};
        ]
      ]
      [.\node[regular]{$16$};
        [.\node[leaf]{$\phantom{00}$};
        ]
        [.\node[leaf]{$\phantom{00}$};
        ]
      ]
    ]
  ]
  [.\node[regular]{$21$};
    [.\node[regular]{$22$};
      [.\node[regular]{$23$};
        [.\node[leaf]{$\phantom{00}$};
        ]
        [.\node[leaf]{$\phantom{00}$};
        ]
      ]
      [.\node[regular]{$26$};
        [.\node[leaf]{$\phantom{00}$};
        ]
        [.\node[leaf]{$\phantom{00}$};
        ]
      ]
    ]
    [.\node[regular]{$\phantom{00}$};
      [.\node[regular]{$\phantom{00}$};
        [.\node[leaf]{$\phantom{00}$};
        ]
        [.\node[leaf]{$\phantom{00}$};
        ]
      ]
      [.\node[regular]{$\phantom{00}$};
        [.\node[leaf]{$\phantom{00}$};
        ]
        [.\node[leaf]{$\phantom{00}$};
        ]
      ]
    ]
  ]
]
\edge[draw=none]; {} ]
\end{tikzpicture}
}\\\hline
\resizebox{\textwidth*\real{0.5}-0.2in}{!}{
%drawTree(levels = 3, s = 2, n = 15, labels = [27, 5, 6, 7, 10, 11, 12, 13, 16, 19, 20, 21, 22, 23, 26])
\tikzstyle{super}=[rectangle, draw]
\tikzstyle{regular}=[circle, draw]
\tikzstyle{leaf}=[ellipse, draw]
\begin{tikzpicture}[level distance = 1.5cm, grow=down]
\Tree[.{$\phantom{aaa}\iddots$}
[.\node[super]{$19,20$};
  [.\node[super]{$10,11$};
    [.\node[super]{$\phantom{0}5,6\phantom{0}$};
      [.\node[leaf]{$27$};
      ]
      [.\node[leaf]{$7\phantom{0}$};
      ]
    ]
    [.\node[regular]{$12$};
      [.\node[leaf]{$13$};
      ]
      [.\node[leaf]{$16$};
      ]
    ]
  ]
  [.\node[regular]{$21$};
    [.\node[regular]{$22$};
      [.\node[leaf]{$23$};
      ]
      [.\node[leaf]{$26$};
      ]
    ]
    [.\node[regular]{$\phantom{00}$};
      [.\node[leaf]{$\phantom{00}$};
      ]
      [.\node[leaf]{$\phantom{00}$};
      ]
    ]
  ]
]
\edge[draw=none]; {} ]
\end{tikzpicture}
} & 
\resizebox{\textwidth*\real{0.5}-0.2in}{!}{
%drawTree(levels = 3, s = 2, n = 15)
\tikzstyle{super}=[rectangle, draw]
\tikzstyle{regular}=[circle, draw]
\tikzstyle{leaf}=[ellipse, draw]
\begin{tikzpicture}[level distance = 1.5cm, grow=down]
\Tree[.{$\phantom{aaa}\iddots$}
[.\node[super]{$10,11$};
  [.\node[super]{$\phantom{0}5,6\phantom{0}$};
    [.\node[super]{$\phantom{0}2,3\phantom{0}$};
      [.\node[leaf]{$1\phantom{0}$};
      ]
      [.\node[leaf]{$4\phantom{0}$};
      ]
    ]
    [.\node[regular]{$7\phantom{0}$};
      [.\node[leaf]{$8\phantom{0}$};
      ]
      [.\node[leaf]{$9\phantom{0}$};
      ]
    ]
  ]
  [.\node[regular]{$12$};
    [.\node[regular]{$13$};
      [.\node[leaf]{$14$};
      ]
      [.\node[leaf]{$15$};
      ]
    ]
    [.\node[regular]{$\phantom{00}$};
      [.\node[leaf]{$\phantom{00}$};
      ]
      [.\node[leaf]{$\phantom{00}$};
      ]
    ]
  ]
]
\edge[draw=none]; {} ]
\end{tikzpicture}
}\\
\end{tabular}
\end{center}
\caption{$T^{(2)}\p{27}$: Prior to pruning (upper left), after Deletion (upper right), after Correction (center left), after Lifting (center right), after Restructuring (lower left), after Relabeling (lower right)}
\label{fig:condel}
\end{figure}
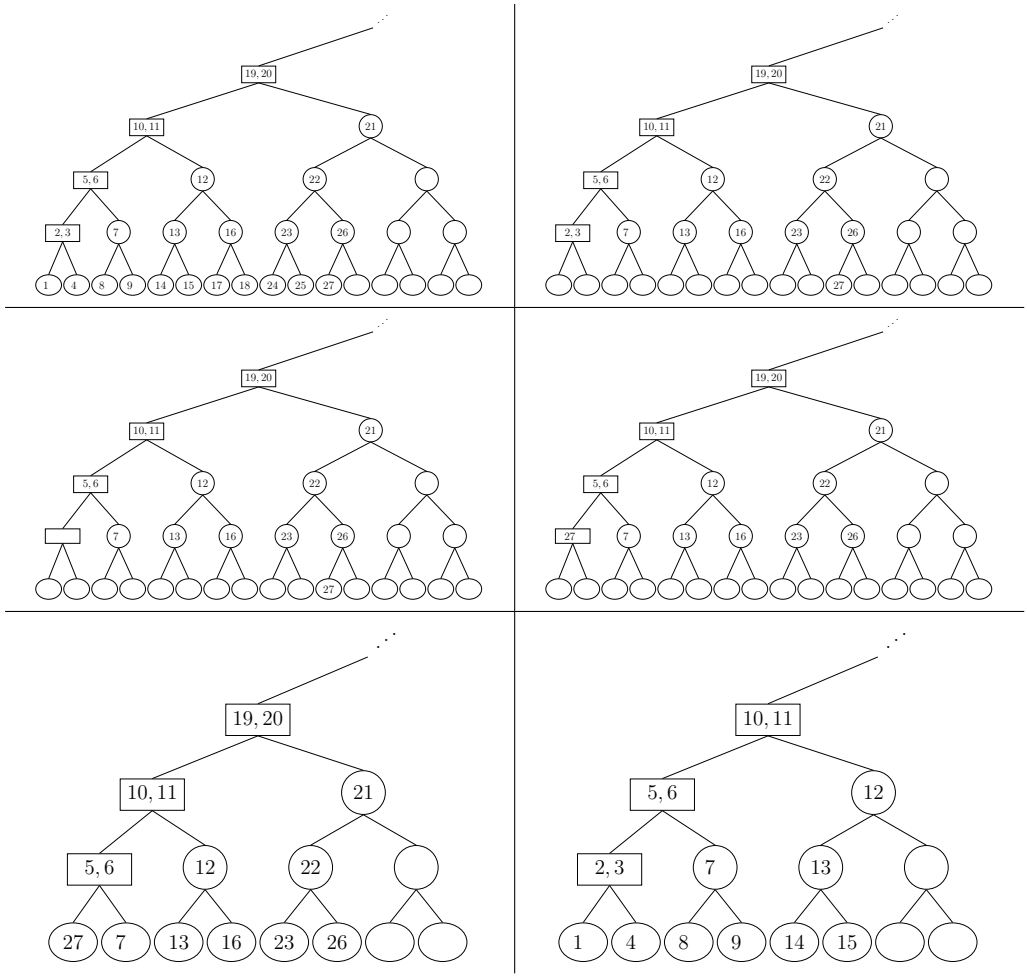
%\begin{figure}
%\begin{center}
%\includegraphics[width=400pt]{contree27rel.eps}
%\end{center}
%\caption{$T_2\p{27}$, after pruning.  Also known as $T_2\p{15}$.}
%\label{fig:conrel}
%\end{figure}
Note that the above process transforms $T^{(s)}\p{n}$ into $T^{(s)}\p{n'}$ for some $n'<n$. 
The Deletion Step deletes $L_{T,s}\p{n-1}$ labels, since deleting all leaf labels that are less than $n$ is equivalent to deleting all labels that are at most $n-1$.  Then, the Correction Step deletes $s$ labels.  The last three steps do not change the number of labels, but they ensure that the result is a valid labeled tree.  So, pruning $T^{(s)}\p{n}$ results in $T^{(s)}\p{n-s-L_{T,s}\p{n-1}}$.  For example, in Figure~\ref{fig:condel}, pruning $T^{(2)}\p{27}$ results in $T^{(2)}\p{15}$.

Furthermore, note that the labeled leaves of the pruned tree $T^{(s)}\p{n-s-L_{T,s}\p{n-1}}$ are precisely the parents of labeled leaves of $T^{(s)}\p{n}$.  The only situation in which a node on level~$1$ would have a label but none of its children would is if it had the label $n$.  But, then its label would have been moved in the Lifting Step.  Since the first child of each regular node to receive a label is its left child, there is a one-to-one correspondence between \typel{0} labeled leaves in $T^{(s)}\p{n}$ and labeled leaves in $T^{(s)}\p{n-s-L_{T,s}\p{n-1}}$.  So, $L_{T,s,0}\p{n}=L_{T,s}\p{n-s-L_{T,s}\p{n-1}}$, as required.  As for \typel{1} leaves, there is the possibility of having a labeled leaf in $T^{(s)}\p{n-s-L_{T,s}\p{n-1}}$ whose \typel{1} child in $T^{(s)}\p{n}$ was unlabeled.  But, this can only happen if the \typel{0} child in $T^{(s)}\p{n}$ was labeled $n$.  Instead, we find that there is a bijection between the labeled \typel{1} leaves of $T^{(s)}\p{n}$ and the labeled leaves of $T^{(s)}\p{n-s-1-L_{T,s}\p{n-2}}$.  If a \typel{1} leaf in $T^{(s)}\p{n}$ has a label, then its \typel{0} sibling has the label one less.  So, to count labeled \typel{1} leaves in $T^{(s)}\p{n}$, we can instead count \typel{0} labeled leaves in $T^{(s)}\p{n-1}$ using the aforementioned bijection.  This results in $L_{T,s,1}\p{n}=L_{T,s}\p{n-s-1-L_{T,s}\p{n-2}}$, as required.
\end{proof}

In the cases where $s<0$, Theorem~\ref{thm:contree} is still true when initial conditions are given by sufficiently many terms of $L_{T,s}\p{n}$. The proof is nearly identical, except that $-s$ \quot{dummy labels} must be added to the regular nodes of level~1 prior to the relabeling step. (The relabeling step replaces these with actual labels.) This ensures that these nodes have labels once they become leaves~\cites{isgur2,ruskeynegative}.

\section{Slow Sequences Connected with Linear Recurrences}\label{sec:treerec}

\subsection{Notation and Terminology}\label{ss:not}
In this section, we consider a particular sort of sequence defined by a linear recurrence.  We consider sequences satisfying recurrences of the form
\[
a_n=\sum_{i=1}^\infty\lambda_ia_{n-i}
\]
subject to the following restrictions:
\begin{itemize}
\item Each $\lambda_i$ is a nonnegative integer, and $\lambda_1$ is a positive integer.
\item For some integer $k\geq1$, for all $i>k$, $\lambda_i=0$.  This value $k$ is the \emph{order} of the recurrence, and is denoted by $k$ going forward.  (We take our sum to $\infty$, rather than just to $k$, as it will sometimes be convenient to refer to $\lambda_i$ for $i>k$.)
\item If $k=1$, then $\lambda_1\geq2$. (We do not allow the recurrence $a_n=a_{n-1}$.)
\item The sequence is generated by the initial conditions $a_i=1$ for all $i\leq0$.
\end{itemize}
Going forward, we denote such a sequence by $\bk{\lambda_1,\lambda_2,\ldots,\lambda_k}$.
Our sequences are similar to Positive Linear Recurrence Sequences (PLRS), differing only in the initial conditions~\cite{millerwangplrs}. 
%Throughout this section, we will consider a generic linear recurrence
%\[
%a_n=\sum_{i=1}^k\lambda_ia_{n-i},
%\]
%where $k$ is the order of the recurrence.  Furthermore, we require the coefficients $\lambda_i$ to be nonnegative integers, and we require $\lambda_1$ to be positive.  Additionally, if we refer to some $\lambda_i$ with $i>k$, its value is taken to be zero.

We now define some auxiliary quantities.  First for each integer $j\geq0$, define
\[
\Lambda_j=\sum_{i=1}^j\lambda_i,
\]
(so $\Lambda_0=0$) and define $\Lambda=\Lambda_k$.  Then, for each integer $0\leq\ell<\Lambda$, define $\mu_\ell$ as the minimum index $j$ where $\Lambda_j>\ell$.  Also, for each nonnegative integer $i$, define $\Delta_i=a_{i+1}-a_i$.

Finally, we define an operator $R$ that allows us to concisely construct nested recurrences of arbitrary depths.  The operator has five parameters:
\begin{description}
\item[$C$:] A symbol representing a nested recurrence relation.
\item[$n$:] A symbol representing the variable in the recurrence relation.
\item[$d$:] A nonnegative integer representing the nesting depth of the resulting recursion.
\item[$s$:] An nonnegative integer that we call a \quot{shift.} 
%An integer, which we will later refer to as a \quot{shift.}  For the rest of Section~\ref{sec:treerec}, $s$ should be treated as nonnegative.  
%Section~\ref{sec:neg} expands to the cases where $s<0$. TODO no
\item[$\ell$:] A nonnegative integer that we call an \quot{offset.}
\end{description}
Define
\[
R\p{C, n, d, s, \ell}=
\begin{cases}
n-\ell & d=0\\
R\p{C, n, d-1, s, \ell}-s-C\p{R\p{C, n, d-1, s, \ell}-1} & d>0.
\end{cases}
\]
%The four parameters 
%In what follows, $L$ denotes a symbol representing a nested recurrence relation, $d$ denotes a nonnegative integer, $s$ denotes any integer (but it should be treated as nonnegative until Section~\ref{sec:neg}), and $j$ denotes a nonnegative integer.
Using the $R$ operator, we can write the recurrence for $C_s$ as
\[
C_s\p{n}=C_s\p{R\p{C_s,n,1,s,0}}+C_s\p{R\p{C_s,n,1,s,1}}.
\]

Our main theorem, Theorem~\ref{thm:main}, is stated below. Theorem~\ref{thm:main} refers to a tree object $T$ that generalizes the Conolly example from Subsection~\ref{ss:tree}. This object $T$ is defined over the next several pages.

\setcounter{temthma}{\value{theorem}}
\begin{theorem}\label{thm:main}
\mainthm
\end{theorem}

\subsection{Trees Defined by Nested Recurrences}\label{ss:treerec}

In Subsection~\ref{ss:tree}, we used a perfect binary tree skeleton in our running example of the tree methodology for finding slow solutions to nested recurrences.  The perfect binary tree skeleton has three properties that make is useful for this purpose:
\begin{description}
\item[Perfection:] All of its leaves are at the same level.
\item[Leaf-Recursivity:] Its structure is preserved by deleting its leaves.
\item[Root-Recursivity:] Each subtree is isomorphic to some left subtree.  (In this case, all subtrees of a given height are isomorphic to one another.)
\end{description}
In our analysis of the Conolly sequence, leaf-recursivity is necessary to make the pruning operation well-defined.  Perfection is used to allow us to discuss levels in the first place. Root-recursivity, while a natural property when thinking about recursive tree constructions, is not actually used in the proof of Theorem~\ref{thm:contree}.  %  TODO counting right children.
With this in mind, we now construct more general trees that can still be used to analyze nested recurrences.  The new trees, which generalize the skeleton for the generalized Conolly recurrences, are perfect and leaf-recursive. They are not root-recursive, though, necessitating a somewhat less intuitive description of their construction.

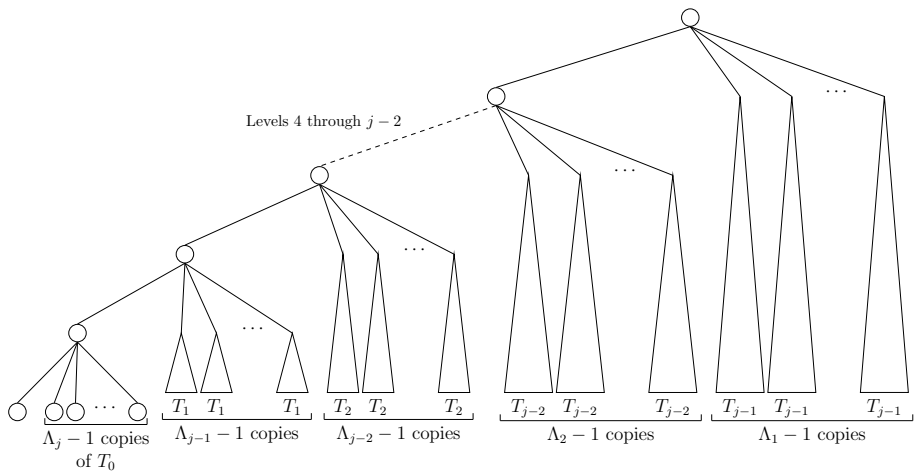
\begin{figure}
\begin{center}
\resizebox{\textwidth-0.5in}{!}{
%\documentclass[border=10pt]{standalone}
%\usepackage{tikz}
%\usepackage{tikz-qtree}
%\usetikzlibrary{calc}
%\begin{document}
\tikzstyle{super}=[rectangle, draw]
\tikzstyle{regular}=[circle, draw]
%\begin{tikzpicture}[level distance = 3cm, grow=down]
%\Tree[.\node[regular]{};
%    [.\node[regular]{};]
%    [.a \edge[roof]; {test}]]
%
%
%\end{tikzpicture}
\begin{tikzpicture}[level distance=50pt]
\begin{scope}[frontier/.style={distance from root=250pt}]
\Tree [.\node[regular]{};
    [.\node[regular]{};
        \edge[dashed] node[auto=right]{{\footnotesize Levels $4$ through $j-2$}};
            [.\node[regular]{}; 
                [.\node[regular]{}; 
                     [.\node[regular]{};
                    \node[regular]{}; \edge[draw=none]; {} \node(t0a)[regular]{}; \node[regular]{}; \edge[draw=none]; {$\cdots$} \node(t0b)[regular]{}; ]
                    \edge[draw=none]; {}
                    [ \edge[roof]; \node(t1a){$T_{1}$}; ]
                 [ \edge[roof]; {$T_{1}$} ]
                  \edge[draw=none]; [.{$\cdots$} ]
                 [ \edge[roof]; \node(t1b){$T_{1}$}; ] ]
                 \edge[draw=none]; {}
                 [ \edge[roof]; \node(t2a){$T_{2}$}; ]
                 [ \edge[roof]; {$T_{2}$} ]
                  \edge[draw=none]; [.{$\cdots$} ]
                 [ \edge[roof]; \node(t2b){$T_{2}$}; ] ]
         \edge[draw=none]; {}     
          \edge[draw=none]; {}         
        [ \edge[roof]; \node(tm2a){$T_{j-2}$}; ]
        [ \edge[roof]; {$T_{j-2}$} ]
         \edge[draw=none]; [.{$\cdots$} ]
        [ \edge[roof]; \node(tm2b){$T_{j-2}$}; ]  ]
   \edge[draw=none]; {}
    [ \edge[roof]; \node(tm1a){$T_{j-1}$}; ] 
    [ \edge[roof]; {$T_{j-1}$} ] 
     \edge[draw=none]; [.{$\cdots$} ]
    [ \edge[roof]; \node(tm1b){$T_{j-1}$}; ] ] 
\end{scope}
%T_0 copies
\path let \p1 = (t0a) in coordinate (t0aa) at (\x1-6,\y1-8);
\path let \p1 = (t0a) in coordinate (t0ab) at (\x1-6,\y1-4);
\path let \p1 = (t0b) in coordinate (t0ba) at (\x1+6,\y1-8);
\path let \p1 = (t0b) in coordinate (t0bb) at (\x1+6,\y1-4);
\draw (t0aa) -- (t0ba)  node[align=center,pos=0.5,below] {$\Lambda_j-1$ copies\\of $T_0$};
\draw (t0aa) -- (t0ab);
\draw (t0ba) -- (t0bb);
%T_1 copies
\path let \p1 = (t1a) in coordinate (t1aa) at (\x1-12,\y1-8);
\path let \p1 = (t1a) in coordinate (t1ab) at (\x1-12,\y1-4);
\path let \p1 = (t1b) in coordinate (t1ba) at (\x1+12,\y1-8);
\path let \p1 = (t1b) in coordinate (t1bb) at (\x1+12,\y1-4);
\draw (t1aa) -- (t1ba)  node[pos=0.5,below] {$\Lambda_{j-1}-1$ copies};
\draw (t1aa) -- (t1ab);
\draw (t1ba) -- (t1bb);
%T_2 copies
\path let \p1 = (t2a) in coordinate (t2aa) at (\x1-12,\y1-8);
\path let \p1 = (t2a) in coordinate (t2ab) at (\x1-12,\y1-4);
\path let \p1 = (t2b) in coordinate (t2ba) at (\x1+12,\y1-8);
\path let \p1 = (t2b) in coordinate (t2bb) at (\x1+12,\y1-4);
\draw (t2aa) -- (t2ba)  node[pos=0.5,below] {$\Lambda_{j-2}-1$ copies};
\draw (t2aa) -- (t2ab);
\draw (t2ba) -- (t2bb);
%T_{j-2} copies
\path let \p1 = (tm2a) in coordinate (tm2aa) at (\x1-18,\y1-8);
\path let \p1 = (tm2a) in coordinate (tm2ab) at (\x1-18,\y1-4);
\path let \p1 = (tm2b) in coordinate (tm2ba) at (\x1+18,\y1-8);
\path let \p1 = (tm2b) in coordinate (tm2bb) at (\x1+18,\y1-4);
\draw (tm2aa) -- (tm2ba)  node[pos=0.5,below] {$\Lambda_2-1$ copies};
\draw (tm2aa) -- (tm2ab);
\draw (tm2ba) -- (tm2bb);
%T_{j-1} copies
\path let \p1 = (tm1a) in coordinate (tm1aa) at (\x1-18,\y1-8);
\path let \p1 = (tm1a) in coordinate (tm1ab) at (\x1-18,\y1-4);
\path let \p1 = (tm1b) in coordinate (tm1ba) at (\x1+18,\y1-8);
\path let \p1 = (tm1b) in coordinate (tm1bb) at (\x1+18,\y1-4);
\draw (tm1aa) -- (tm1ba)  node[pos=0.5,below] {$\Lambda_1-1$ copies};
\draw (tm1aa) -- (tm1ab);
\draw (tm1ba) -- (tm1bb);
\end{tikzpicture}
%\end{document}
}
\end{center}
\caption{A generic tree $T_j$}
\label{fig:gentreei}
\end{figure}

%TODO subscript on T currently means Tj and label
Let $\bk{\lambda_1,\lambda_2,\ldots,\lambda_k}$ be a sequence as described in Subsection~\ref{ss:not}.  We will now define a sequence of finite trees $T_0\bk{\lambda_1,\lambda_2,\ldots,\lambda_k}, T_1\bk{\lambda_1,\lambda_2,\ldots,\lambda_k},\ldots$ and a skeleton $T\bk{\lambda_1,\lambda_2,\ldots,\lambda_k}$.  (Going forward, we drop the bracketed part of the notation if the recurrence in question is generic or otherwise clear.)
\begin{defin}\label{def:ti}
Define a sequence $T_0,T_1,\ldots$ of trees as follows:
\begin{itemize}
\item $T_0$ consists of a single node.
\item $T_j$ has $j+1$ levels, numbered from $0$ to $j$ from bottom to top.  The leftmost node on each level is said to be \emph{special}.
%\item For every level $i>0$, the special node at that level has a left child, namely special child at the level below it.
\item The special node on level $i>0$ has $\Lambda_{j-i+1}$ children.  The leftmost of these children is the special node on level $i-1$.  The rest of these children are roots of copies of $T_{i-1}$.
%\item $T_0$ consists of a single node.
%\item For $j>0$, $T_j$ has $j+1$ levels (numbered from $0$ to $j$ from bottom to top).  The leftmost node at level $i>0$ has $\Lambda_{j-i+1}-1$ copies of $T_{i-1}$ as subtrees, all to the right of its leftmost child.
\end{itemize}
\end{defin}
See Figure~\ref{fig:gentreei} for a generic diagram of $T_j$.  Note that these trees are not skeletons and have only one type of node, always depicted as a circle.

\begin{figure}
\begin{center}
%TODO https://tex.stackexchange.com/questions/208771/how-to-draw-triangles-as-subtree-with-the-forest-package
%\includegraphics[width=300pt]{hof10000.eps}
\resizebox{\textwidth-0.5in}{!}{
%\documentclass[border=10pt]{standalone}
%\usepackage{tikz}
%\usepackage{tikz-qtree}
%\usetikzlibrary{calc}
%\begin{document}
\tikzstyle{super}=[rectangle, draw]
\tikzstyle{regular}=[circle, draw]
%\begin{tikzpicture}[level distance = 3cm, grow=down]
%\Tree[.\node[regular]{};
%    [.\node[regular]{};]
%    [.a \edge[roof]; {test}]]
%
%
%\end{tikzpicture}
\begin{tikzpicture}[level distance=50pt]
\begin{scope}[frontier/.style={distance from root=300pt}]
\Tree [.{$\phantom{aaa}\iddots$} [.\node[super]{\phantom{$0000$}};
    [.\node[super]{\phantom{$0000$}};
        \edge[dashed] node[auto=right]{{\footnotesize Levels $4$ through $j-2$}};
            [.\node[super]{\phantom{$0000$}}; 
                [.\node[super]{\phantom{$0000$}}; 
                     [.\node[super]{\phantom{$0000$}};
                    \node[regular]{}; \edge[draw=none]; {} \node(t0a)[regular]{}; \node[regular]{}; \edge[draw=none]; {$\cdots$} \node(t0b)[regular]{}; ]
                    \edge[draw=none]; {}
                    [ \edge[roof]; \node(t1a){$T_{1}$}; ]
                 [ \edge[roof]; {$T_{1}$} ]
                  \edge[draw=none]; [.{$\cdots$} ]
                 [ \edge[roof]; \node(t1b){$T_{1}$}; ] ]
                 \edge[draw=none]; {}
                 [ \edge[roof]; \node(t2a){$T_{2}$}; ]
                 [ \edge[roof]; {$T_{2}$} ]
                  \edge[draw=none]; [.{$\cdots$} ]
                 [ \edge[roof]; \node(t2b){$T_{2}$}; ] ]
         \edge[draw=none]; {}     
          \edge[draw=none]; {}         
        [ \edge[roof]; \node(tm2a){$T_{j-2}$}; ]
        [ \edge[roof]; {$T_{j-2}$} ]
         \edge[draw=none]; [.{$\cdots$} ]
        [ \edge[roof]; \node(tm2b){$T_{j-2}$}; ]  ]
   \edge[draw=none]; {}
    [ \edge[roof]; \node(tm1a){$T_{j-1}$}; ] 
    [ \edge[roof]; {$T_{j-1}$} ] 
     \edge[draw=none]; [.{$\cdots$} ]
    [ \edge[roof]; \node(tm1b){$T_{j-1}$}; ] ] 
    \edge[draw=none]; {} ]
\end{scope}
%T_0 copies
\path let \p1 = (t0a) in coordinate (t0aa) at (\x1-6,\y1-8);
\path let \p1 = (t0a) in coordinate (t0ab) at (\x1-6,\y1-4);
\path let \p1 = (t0b) in coordinate (t0ba) at (\x1+6,\y1-8);
\path let \p1 = (t0b) in coordinate (t0bb) at (\x1+6,\y1-4);
\draw (t0aa) -- (t0ba)  node[align=center,pos=0.5,below] {$\Lambda-1$ copies\\of $T_0$};
\draw (t0aa) -- (t0ab);
\draw (t0ba) -- (t0bb);
%T_1 copies
\path let \p1 = (t1a) in coordinate (t1aa) at (\x1-12,\y1-8);
\path let \p1 = (t1a) in coordinate (t1ab) at (\x1-12,\y1-4);
\path let \p1 = (t1b) in coordinate (t1ba) at (\x1+12,\y1-8);
\path let \p1 = (t1b) in coordinate (t1bb) at (\x1+12,\y1-4);
\draw (t1aa) -- (t1ba)  node[pos=0.5,below] {$\Lambda-1$ copies};
\draw (t1aa) -- (t1ab);
\draw (t1ba) -- (t1bb);
%T_2 copies
\path let \p1 = (t2a) in coordinate (t2aa) at (\x1-12,\y1-8);
\path let \p1 = (t2a) in coordinate (t2ab) at (\x1-12,\y1-4);
\path let \p1 = (t2b) in coordinate (t2ba) at (\x1+12,\y1-8);
\path let \p1 = (t2b) in coordinate (t2bb) at (\x1+12,\y1-4);
\draw (t2aa) -- (t2ba)  node[pos=0.5,below] {$\Lambda-1$ copies};
\draw (t2aa) -- (t2ab);
\draw (t2ba) -- (t2bb);
%T_{j-2} copies
\path let \p1 = (tm2a) in coordinate (tm2aa) at (\x1-18,\y1-8);
\path let \p1 = (tm2a) in coordinate (tm2ab) at (\x1-18,\y1-4);
\path let \p1 = (tm2b) in coordinate (tm2ba) at (\x1+18,\y1-8);
\path let \p1 = (tm2b) in coordinate (tm2bb) at (\x1+18,\y1-4);
\draw (tm2aa) -- (tm2ba)  node[pos=0.5,below] {$\Lambda-1$ copies};
\draw (tm2aa) -- (tm2ab);
\draw (tm2ba) -- (tm2bb);
%T_{j-1} copies
\path let \p1 = (tm1a) in coordinate (tm1aa) at (\x1-18,\y1-8);
\path let \p1 = (tm1a) in coordinate (tm1ab) at (\x1-18,\y1-4);
\path let \p1 = (tm1b) in coordinate (tm1ba) at (\x1+18,\y1-8);
\path let \p1 = (tm1b) in coordinate (tm1bb) at (\x1+18,\y1-4);
\draw (tm1aa) -- (tm1ba)  node[pos=0.5,below] {$\Lambda-1$ copies};
\draw (tm1aa) -- (tm1ab);
\draw (tm1ba) -- (tm1bb);
\end{tikzpicture}
%\end{document}
}
\end{center}
\caption{Levels $0$ through $j$ of the generic skeleton $T$}
\label{fig:gentree}
\end{figure}
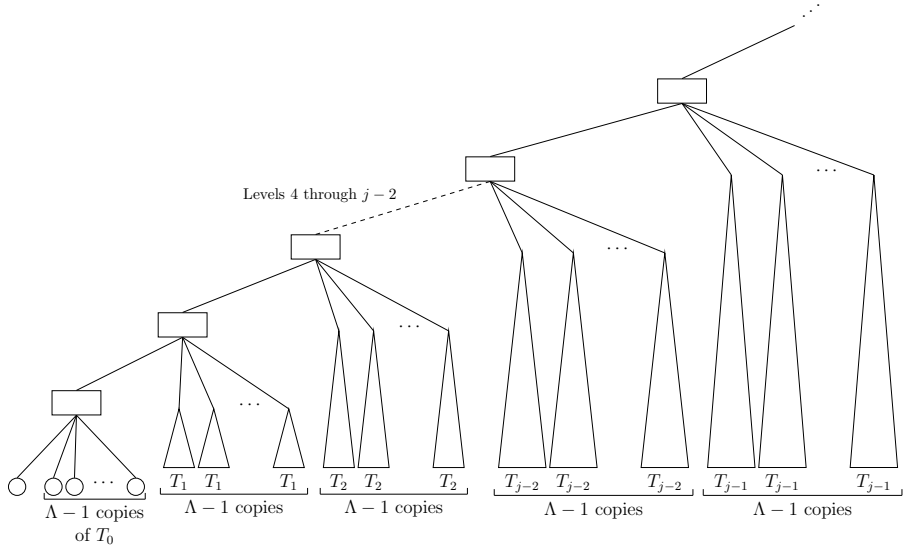

\begin{defin}\label{def:t}
Define a skeleton $T$ as follows:
\begin{itemize}
\item The node at the lower left is a leaf.
\item The leftmost node at every level other than the bottom is a supernode.
\item The $i^{th}$ supernode has $\Lambda-1$ copies of $T_{i-1}$ as subtrees, all to the right of its leftmost child (the special node below it).
\end{itemize}
\end{defin}
See Figure~\ref{fig:gentree} for a generic diagram of $T$.
Note that $T\bk{2}$ is the perfect binary tree example from Subsection~\ref{ss:tree}.  Also, observe that, by construction, every node has at most $\Lambda$ children.

%TODO make sure labeling and $L_{T,s}\p{n}$ are defined earlier

%TODO distinction between skeleton w/supernodes and without?% TODO define index of a supernode (1st, 2nd, etc)

We have the following claims about these trees:
\begin{prop}\label{prop:tperf}
The trees $T_j$ and $T$ are perfect.
\end{prop}
\begin{proof}
First, we show inductively that $T_j$ is perfect.  Since $T_0$ is a single vertex, it is perfect.  Now, suppose each $T_m$ is perfect for $m<j$.  The leftmost leaf of $T_j$ is at level~$0$, and the root is at level~$j$.  For each other leaf in $T_j$, define its \emph{pedigree} to be the lowest level for which it is descended from the special node at that level.  Since all leaves descend from the root, which is special, each leaf has a pedigree.  A leaf of pedigree $i$ is a leaf in a copy of $T_{i-1}$.  By induction, every leaf of $T_{i-1}$ is $i-1$ levels below its root.  This means that such a leaf in $T_j$ is $i$ levels below its pedigree-defining node.  But, its pedigree-defining node is on level $i$, so a leaf of pedigree $i$ is on level~$0$, as required.%Every leaf below the leftmost node on level $i$ is a leaf of a copy of $T_{i-1}$.  So, by induction, every such leaf is $i$ levels below level $i$, or on level~$0$, as required.

Finally, we show that every leaf of $T$ is on level~$0$.  The leftmost leaf is on level~$0$ by definition.  Every leaf below the $i^{th}$ supernode is a leaf of a copy of $T_{i-1}$.  So, every such leaf is $i$ levels below level $i$, or on level~$0$, as required.
\end{proof}

\begin{prop}\label{prop:tilr}
For $j>0$, deleting the leaves from $T_j$ results in $T_{j-1}$.
\end{prop}
\begin{proof}
The proof is by induction on $j$.  Note that $T_1$ has two levels, so deleting its leaves results in a single vertex, which is $T_0$.  Now, suppose that, for $m<j$, deleting the leaves from $T_m$ results in $T_{m-1}$.  Consider what happens when we delete the leaves from $T_j$.  The special node at level $1$ in $T_j$ has all of its children removed.  The special node at level $i>1$ in $T_j$ has $\Lambda_{j-i+1}$ copies of $T_{i-1}$ below it, to the right of its special child.  Deleting the leaves from $T_j$ deletes the leaves from these subtrees, so a node having $\Lambda_{j-i+1}$ copies of $T_{i-1}$ as subtrees in $T_j$ has $\Lambda_{j-i+1}$ copies of $T_{i-2}$ as subtrees once the leaves are deleted from $T_j$.  Reindexing the levels so that the new bottom level is level $0$ completes the proof.
\end{proof}

\begin{prop}\label{prop:tlr}
The skeleton $T$ is leaf-recursive.
\end{prop}
\begin{proof}
Consider what happens when we delete the leaves from $T$.  The first supernode has all of its children removed.  For $i>1$, the $i^{th}$ supernode has $\Lambda-1$ copies of $T_{i-1}$ as subtrees, to the right of its supernode child.  Deleting the leaves from $T$ deletes the leaves from these subtrees.  So, by Proposition~\ref{prop:tilr}, a node having $\Lambda-1$ copies of $T_{i-1}$ as subtrees in $T$ has $\Lambda-1$ copies of $T_{i-2}$ as subtrees once the leaves are deleted from $T$.  Reindexing the levels so that the new bottom level is level $0$ completes the proof.
\end{proof}

%We will frequently refer to nodes based on how far to the left they are relative to their siblings.  We will call a node 
Recall (Definition~\ref{def:type}) that a node is called \typel{\ell} if it has exactly $\ell$ siblings to its left.  %Our term \emph{special} from Definition~\ref{def:ti} will henceforth be called \emph{\typel{0}}.  TODO types earlier TODO this is wrong
We have the following results about types and child counts:
\begin{prop}\label{prop:tchildcount}
In $T$:
\begin{enumerate}
\item Every supernode has $\Lambda$ children.\label{it:super}
\item Every non-supernode has $\Lambda_{j+1}$ children, where $j$ is the distance to the most recent non-\typel{0} ancestor.  (If the node in question is not \typel{0}, $j=0$.)\label{it:nosuper}
\end{enumerate}
\end{prop}
\begin{proof}
Item~\ref{it:super} is clear from the construction.  For Item~\ref{it:nosuper}, observe that any non-\typel{0} node is the root of a copy of some $T_m$.  By the construction of $T_m$, a node $j$ levels to the left below the root of $T_m$ has $\Lambda_{j+1}$ children, as required.
\end{proof}

\section{Main Result}

Now that we have defined our skeleton, we restate and prove our main result. The proof generalizes the proof of Theorem~\ref{thm:contree}.

\setcounter{temthm}{\value{theorem}}
\setcounter{theorem}{\value{temthma}}
\begin{theorem}
\mainthm
\end{theorem}
\setcounter{theorem}{\value{temthm}}

\begin{proof}
%(Here, we prove the case $s\geq0$.  The case $s<0$ is left to Section~\ref{sec:neg}.)
%
We argue by induction on $n$, with the base case provided by the initial condition.  Suppose $n>a_k$, and suppose that $C_a\p{\ninduct}=L_{T,s}\p{\ninduct}$ for all $\ninduct<n$.

The value of $L_{T,s}\p{n}$ is defined as the number of labeled leaves in $T^{(s)}\p{n}$.  For each $0\leq\ell<\Lambda$, define $L_{T,s,\ell}\p{n}$ as the number of labeled \typel{\ell} leaves in $T^{(s)}\p{n}$.  Clearly,
\[
L_{T,s}\p{n}=\sum_{\ell=0}^{\Lambda-1}L_{T,s,\ell}\p{n}.
\]

We now claim that $L_{T,s,\ell}\p{n}=C_a\p{R\p{C_a, n,\mu_\ell, s, \ell}}$.  We define a pruning process for $T^{(s)}\p{n}$ analogously to the pruning process in Subsection~\ref{ss:tree}:
\begin{description}
\item[Deletion Step:] Delete all leaf labels in $T^{(s)}\p{n}$ that are less than $n$.
\item[Correction Step:] Delete all $s$ labels from the first supernode.
\item[Lifting Step:] Move the label $n$ to the first supernode.
\item[Relabeling Step:] The old first supernode becomes a leaf, and relabel the currently labeled nodes in preorder.
\end{description}
The Deletion Step deletes $L_{T,s}\p{n-1}$ labels, the Correction Step deletes $s$ labels, and the last two steps do not change the number of labels but ensure that the result is a valid labeled tree.  So, pruning $T^{(s)}\p{n}$ results in $T^{(s)}\p{n-L_{T,s}\p{n-1}-s}$.

As a consequence of Proposition~\ref{prop:tchildcount}, the $\mu_\ell-1$ immediate ancestors of every \typel{\ell} leaf are \typel{0}.  Also, every node at the end of a chain of $\mu_\ell-1$ consecutive \typel{0} nodes has a \typel{\ell} child.  Since every non-leaf node has a \typel{0} child, we have 
a bijection between labeled nodes on Level-$\mu_\ell$ and the \typel{\ell} leaves in $T^{(s)}$ that are either themselves labeled at most $n$ or have a sibling labeled at most $n$. 
%Note that it is not a bijection between the labeled nodes on Level-$\mu_\ell$ and the \typel{\ell} leaves labeled at most $n$
%that the number of \typel{\ell} leaves that are either themselves labeled at most $n$ or have a sibling labeled at most $n$ is the same as the number of leaves present when $T^{(s)}\p{n-\ell}$ is pruned $\mu_\ell$ times.  This number is almost $L_{T,s,\ell}\p{n}$, but it is one too large if the leaf labeled $n$ has a \typel{\ell} sibling and is itself of type less than $\ell$.  To avoid this issue, we can, instead of pruning $T^{(s)}\p{n}$, prune $T^{(s)}\p{n-\ell}$.  This operation excludes the aforementioned ineligible siblings by pruning a tree without any leaves in the problematic set of siblings, and it does not throw away any leaves we care about in other cases. TODO justify.
We are only interested in the \typel{\ell} leaves that themselves are labeled at most $n$. So, we can look instead look at \typel{0} leaves in $T^{(s)}$ with label at most $n-\ell$ that have a \typel{\ell} sibling.  This gives the more useful bijection between these nodes and their ancestors on Level-$\mu_\ell$. (At most one node on Level-$\mu_\ell$ is excluded here that was included before.)

From here, we obtain that $L_{T,s,\ell}\p{n}$ equals the number of leaves in the tree resulting from pruning $T^{(s)}\p{n-\ell}$ exactly $\mu_\ell$ times.  (Each pruning after the first matches \typel{0} nodes with their parents, so we need not make further adjustments to later prunings.) We now claim that the number of nodes in the tree resulting from pruning $T_{n-\ell,s}$ exactly $d$ times is $R\p{C_a, n, d, s, \ell}$.  We argue by induction on $d$.  If $d=1$, we prune once, which results in $T^{(s)}\p{n-L_{T,s}\p{n-1}-s}$.  By our earlier inductive hypothesis, $L_{T,s}\p{n-1}=C_a\p{n-1}$, so the number of nodes is $n-s-C_a\p{n-1}$.  This equals $R\p{C_a, n, 1, s, \ell}$, as required.  Now, suppose $d>1$, and suppose that the tree resulting from pruning $T^{(s)}\p{n-\ell}$ exactly $d-1$ times has $R\p{C_a, n, d-1, s, \ell}$ nodes, meaning it is $T^{(s)}\p{R\p{C_a, n, d-1, s, \ell}}$.  Pruning it again results in $T^{(s)}\p{R\p{C_a, n, d-1, s, \ell}-L_{T,s}\p{R\p{C_a, n, d-1, s, \ell}-1}-s}$.  By our earlier inductive hypothesis, $L_{T,s}\p{R\p{C_a, n, d-1, s, \ell}-1}=C_a\p{R\p{C_a, n, d-1, s, \ell}-1}$, so the number of nodes is $R\p{C_a, n, d-1, s, \ell}-s-C_a\p{R\p{C_a, n, d-1, s, \ell}-1}=R\p{C_a,n,d,s,\ell}$, as required.

Now, by our earlier inductive hypothesis, the number of leaves in the tree resulting from pruning $T^{(s)}\p{n-\ell}$ exactly $\mu_\ell$ times, and hence $L_{T,s,\ell}\p{n}$, equals $C_a\p{R\p{C_a,n,d,s,\ell}}$.  We therefore have
\[
L_{T,s}\p{n}=\sum_{\ell=0}^{\Lambda-1}L_{T,s,\ell}\p{n}
=\sum_{\ell=0}^{\Lambda-1}C_a\p{R\p{C_a,n,d,s,\ell}},
\]
as required.
\end{proof}

As was the case in Subsection~\ref{ss:tree}, Theorem~\ref{thm:main} is still true when $s<0$, provided sufficiently many values of $L_{T,s}\p{n}$ are used as initial conditions. The sequences where $s<0$ lack the combinatorial interpretation in the upcoming section, so we end our discussion of them here.

Three examples of Theorem~\ref{thm:main} being applied can be found in Appendix~\ref{sec:examples}.

\section{Another Combinatorial Interpretation: A Generalization of Zeckendorf Representations}\label{sec:genzeck}

Our sequences have another combinatorial interpretation, which generalizes the following property of the Conolly sequence that we mentioned in the introduction: The number of occurrences of the number $N$ in the Conolly sequence is one plus the number of zeroes at the end of the binary representation of $N$~\cite[A046699]{oeis}. Formally, the \emph{base-$b$} expansion of the nonnegative integer $N$ is the unique sequence of integers $\dseq{d}{m}$ called \emph{digits} with $0\leq d_i<b$ and
\[
\sum_{i=0}^md_ib^i=N.
\]
More generally, given a strictly increasing sequence $\seq{a_k}_{k\geq0}$ of positive integers with $a_0=1$, we can discuss representing $N$ by a sum of the form
\[
\sum_{i=0}^md_ia_i.
\]
Base-$b$ representations are obtained from the case $a_i=b^i$. The most famous non-base-$b$ system of this form is the \emph{Zeckendorf representation}, obtained by taking $a_i=F_{i+2}$, where $F_k$ denotes the $k^{th}$ Fibonacci number. Every positive integer can be represented in this system, and the representation is unique provided that each digit is $0$ or $1$ and provided that no two consecutive digits are $1$. The Zeckendorf representation of $N$ can be obtained by repeatedly greedily subtracting the largest possible Fibonacci number from $N$ while keeping the result nonnegative~\cite{zeck}. A similar greedy algorithm can be used to obtain the base-$b$ representation of $N$. Obtaining representations from other sequences in a similar manner is sometimes known as \emph{digital representation} or \emph{generalized Zeckendorf representation}~\cites{drmota,millerwangplrs}.

In our discussion that follows, we use sequences of digits to represent positive integers, and we refer to this process as a generalization of Zeckendorf representations. This is, in fact, a \emph{different} generalization of Zeckendorf representations than the natural one above. The generalizations coincide for base-$b$ and Zeckendorf representations, but they are different in almost all other cases. (See Proposition~\ref{prop:placeval}.)

%TODO cite some stuff on Zeckendorf and generalized Zeckendorf that isn't this.

%TODO word notation, including weird indexing from end and from 0 with brackets

\begin{defin}\label{def:valid}
Given a recurrence $a=\bk{\lambda_1,\lambda_2,\ldots,\lambda_k}$
and a positive integer $M$, the 
sequence of digits $\dseq{d}{M}$ 
is 
%word $w$
\emph{$a$-valid} if it has the following properties:
\begin{itemize}
%\item $w\in\st{0,1,2,\ldots,\Lambda-1}^*$.
\item For each $i$, $0\leq d_i<\Lambda$. (All digits are nonnegative and are less than $\Lambda$.)
%\item $w$ starts with a symbol other than $0$.  (In particular, the empty word is not $a$-valid.)
\item $d_M\neq0$. (The first digit is nonzero.)
%\item Whenever the symbol $\ell$ appears in $w$, it must either be the first symbol in $w$, or it must be preceded by at least $\mu_\ell-1$ consecutive zeroes.
\item If $d_i=\ell$, then either $i=M$ or $d_{i+1}=d_{i+2}=\cdots=d_{i+\mu_\ell-1}=0$. (The digit $\ell$ is preceded by at least $\mu_\ell-1$ consecutive zeroes.)
\end{itemize}
\end{defin}
%Note that our sequences are written as $d_0,d_1,\ldots,d_M$, but we say that digit $d_M$ is the \emph{first} digit, $d_0$ is the \emph{last} digit, and that digit $d_{i+1}$ \emph{precedes} digit $d_i$. These conventions align with usual conventions when writing base-$b$ expansions of numbers.
%Note that if the recurrence $a$ is obvious from context, we simply call a word \emph{valid}.

Let $Z_a$ denote the set of all $a$-valid 
%words.  Order $Z_a$ so that shorter words appear before longer words, and order words of the same length lexicographically.  Define the \emph{$a$-Zeckendorf} representation of the positive integer $N$ to be the $N^{th}$ word in the ordering on $Z_a$.  TODO not same generalization TODO Zeck for Fib
sequences. Order $Z_a$ so that shorter sequences appear before longer ones, and order sequences of the same length lexicographically. %, reading from $d_M$ to $d_0$. 
Define the \emph{$a$-Zeckendorf} representation of the positive integer $N$ to be the $N^{th}$ sequence in the ordering on $Z_a$. Note that if $a=\bk{b}$, then the $a$-Zeckendorf representation of $N$ is the base-$b$ representation of $N$. Furthermore, if $a=\bk{1,1}$, then the $a$-Zeckendorf representation of $N$ is the Zeckendorf representation of $N$.

We have the following claims about $a$-Zeckendorf representations.

\begin{prop}\label{prop:validz}
Let $m$ be a nonnegative integer.  There are $\Lambda_{m+1}$ 
%symbols that are allowed to appear in a valid string after a run of $m$ consecutive zeroes.
digits that are allowed to appear in a valid sequence after a run of exactly $m$ consecutive zeroes.
\end{prop}
\begin{proof}
%A symbol $\ell$ is allowed to appear after a run of $m$ consecutive zeroes precisely when $\mu_\ell-1\leq m$, i.e.\ when $\mu_\ell\leq m+1$.  The number $\mu_\ell$ is defined as the minimum number $j$ such that $\Lambda_j>\ell$.  So, $\mu_\ell\leq m+1$ precisely when $\ell\leq\Lambda_{m+1}$, as required.
A digit $\ell$ is allowed to appear after a run of $m$ consecutive zeroes precisely when $\mu_\ell-1\leq m$, i.e.\ when $\mu_\ell\leq m+1$.  The number $\mu_\ell$ is defined as the minimum number $j$ such that $\Lambda_j>\ell$.  So, $\mu_\ell\leq m+1$ precisely when $\ell\leq\Lambda_{m+1}$, as required.
\end{proof}

\begin{prop}\label{prop:zcount}
Let $t$ be any integer.  The number of $a$-valid
%words of length at most $t$ is $a_t-1$.
sequences with at most $t$ digits is $a_t-1$.
\end{prop}
\begin{proof}
We proceed by induction on $t$.  Every $a$-valid sequence has at least $1$ digit, so the number of $a$-valid sequences with at most $t$ digits is $0=a_t-1$ whenever $t\leq0$.  Now, suppose $t>0$, and consider an $a$-valid sequence $D=\dseq{d}{M}$ of length at most $t$ with $d_0=\ell$. Furthermore, suppose that the number of $a$-valid sequences of length at most $t'<t$ is $a_{t'}-1$.  Let $D'=\dseqend{d}{M}{2}{1}$.  We must have that either $D'$ is empty or that $D'$ is an $a$-valid sequence with at most $t-1$ digits that ends in at least $\mu_\ell-1$ zeroes.  An $a$-valid sequence with at most $t-1$ digits that ends in at least $\mu_\ell-1$ zeroes corresponds to any $a$-valid sequence with at most $t-\mu_\ell$ digits, followed by $\mu_\ell-1$ zeroes.  So, the number of $a$-valid sequences with at most $t$ digits that have $d_0=\ell$ is one plus the number of $a$-valid sequences with $t-\mu_\ell$ digits, which is $a_{t-\mu_\ell}-1$.  There is one exception to all of this: If $\ell=0$, $D'$ cannot be empty.  Putting everything together yields that the number of $a$-valid sequences with at most $t$ digits:
\begin{equation}\label{eq:tmell}
-1+\sum_{\ell=0}^{\Lambda-1}\pb{1+a_{t-\mu_\ell}-1}=-1+\sum_{\ell=0}^{\Lambda-1}a_{t-\mu_\ell}.
\end{equation}
Since each integer $i\geq1$ appears as $\mu_\ell$ for exactly $\lambda_i$ different values of $\ell$, Equation~\ref{eq:tmell} is the same as
\[
-1+\sum_{i=1}^k\lambda_ia_{t-i}=a_t-1,
\]
as required.
%Each valid string of length at most $t$ can be created 
%in exactly one of the following $\Lambda$ ways:
%\begin{itemize}
%\item Start with an empty string or a valid string of length at most $t-\mu_0$ and concatenate a $0$ to the end.
%\item Start with an empty string or a valid string of length at most $t-\mu_1$ and concatenate a $1$ to the end.
%\end{itemize}
%$t$ is $\Delta_t$.  Observe that the sequence of $\Delta$ values satisfies the same recurrence as the sequence of $a$ values.  So, we will show that the sequence enumerating the number of valid words of length $t$ satisfies the same recurrence as $a$ while having initial conditions $\Delta_0,\Delta_1,\ldots,\Delta_
\end{proof}

\begin{cor}\label{cor:zcount}
Let $N$ be a positive integer, and let $t$ be such that $a_t\leq n<a_{t+1}$.  Then, the $a$-Zeckendorf representation of $N$ has $t+1$ digits.
\end{cor}
\begin{proof}
By Proposition~\ref{prop:zcount}, the number of $a$-valid sequences with at most $t$ digits is $a_t-1$, and the number of $a$-valid sequences with at most $t+1$ digits is $a_{t+1}-1$.  The former number is less than $N$, and the latter number is greater than or equal to $N$.  So, the $N^{th}$ sequence in order must have exactly $t$ digits, as required.
\end{proof}

We now formally state and prove our earlier claim that $a$-Zeckendorf representations and previous notions of digital representations coincide for base-$b$ representations and Zeckendorf representations but differ in almost all other cases. For simplicity, we introduce the following definition.

%In our earlier discussion of Zeckendorf representations, we claimed that $a$-Zeckendorf representations are not the as the more typical sort of generalized Zeckendorf representation. We now explore this distinction.
\begin{defin}\label{def:placeval}
A digital representation system is a \emph{place value system} if there exists a strictly increasing sequence of integers $c_0,c_1,c_2,\ldots$, called \emph{place values}, where $c_0=1$ and where the number represented as %$d_jd_{j-1}\cdots d_2d_1d_0$ is
$\dseq{d}{M}$ is
\[
\sum_{i=0}^Md_ic_i.
\]
\end{defin}
This leads to the following proposition.
\begin{prop}\label{prop:placeval}
The $a$-Zeckendorf representation system is a place value system if and only if one of the following conditions holds:
\begin{itemize}
\item $k=1$
%\item $\lambda_k\geq2$ if $k=1$
\item $k>1$ and $\Lambda=2$
\end{itemize}
If it is a place value system, the place values are the numbers $a_i$.
\end{prop}
\begin{proof}
First, it is clear that if the $a$-Zeckendorf system is a place value system, the place values must be the numbers $a_i$, because the representation of $a_i$ is always a $1$ followed by $i$ zeroes. We now prove the rest of the claim.
\iffpf{
%Forward direction
Here, we prove the contrapositive. Suppose that $k>1$ and that $\Lambda\geq3$. We claim that the number represented as $20$ in the $a$-Zeckendorf system is not $2a_1$, which suffices to show that this is not a place value system. For simplicity of notation, let $L=\Lambda-1$. Since $k>1$, the sequence $1L$ is not $a$-valid. This means that there are fewer than $a_1+\Lambda-1$ $a$-valid sequences preceding $20$ (the $-1$ is to exclude the invalid sequence $0$). But, it is always the case that $a_1=\Lambda$, so there are fewer than $2a_1-1$ $a$-valid sequences preceding $20$. This means that $20$ cannot represent $2a_1$, as required.
}
{
%Reverse direction
First, suppose $k=1$. Whenever $k=1$, we require $\lambda_1\geq2$. For any of these cases, any 
sequence of nonnegative digits that are less than $\lambda_1$ 
%word over the appropriate alphabet 
is $a$-valid. This results in the base-$\lambda_1$ representation system, which is a familiar place value system.

Now, suppose $k>1$ and $\Lambda=2$. Since we always have $\lambda_1\geq1$ and $\lambda_k\geq1$, our recurrence must be $a_n=a_{n-1}+a_{n-k}$.  If $k=2$, this is the Fibonacci sequence, and the $a$-Zeckendorf representation is the Zeckendorf representation. The rest of this proof generalizes the proof of Zeckendorf's Theorem.

Consider the place value system with the numbers $a_i$ as place values. In this system, representation of a given integers is not unique, so consider on the \quot{greedy} representation. To represent the number $N$, use the following algorithm. Start with $M=N$, and continue until $M=0$. Find the index $j$ such that $a_j\leq M<a_{j+1}$. Set $d_j=1$, then replace $M$ by $M-a_j$. 
Since $a_j\leq M<a_{j+1}$, we have $M-a_j<a_{j+1}-a_j$. We know that $a_{j+1}=a_j+a_{j-k+1}$, so $a_{j+1}-a_j=a_{j-k+1}$. This shows that $M-a_j<a_{j-k+1}$. Since both systems use only zeroes and ones, this shows that any two ones are separated by at least $k-1$ zeroes. Therefore, this representation is the same as the $a$-Zeckendorf representation, as required.
%We claim that this representation is the same as the $a$-Zeckendorf representation. To prove this, it will suffice to show that $M-a_j<a_{j-k+1}$. Both systems only use zeroes and ones, and this would show that any two ones are separated by at least $k-1$ zeroes.
}
\end{proof}

\subsection{Connection with Trees}
We now connect $a$-Zeckendorf representations with the skeletons from Subsection~\ref{ss:treerec}.

\begin{prop}\label{prop:ztree}
Let $t$ be a positive integer, and let $T$ be the skeleton for the recurrence $a$ (Definition~\ref{def:t}).  Consider the following pair of maps, one from $a$-valid sequences with at most $t$ digits to non-leftmost leaves below the $t^{th}$ supernode of $T$ and the other from these leaves to these $a$-valid sequences:
\begin{description}
\item[Digits to Leaves:] Let $\dseq{d}{M}$ be an $a$-valid sequence with $M\leq t$. Start with the $M^{th}$ supernode as the current node, and read the digits from $d_M$ to $d_0$.  As symbol $\ell$ is read from the sequence, move from the current node to its \typel{\ell} child and make that the new current node.
\item[Leaves to Digits:] Let $v$ be a leaf below the $t^{th}$ supernode that is not the leftmost leaf.  Start with $v$ as the current node and 
$i=0$. Until the current node is a supernode, let $d_i$ equal the type of the current node, and then replace the current node by its parent.
%$w$ as an empty word.  Until the current node is a supernode, concatenate the type of the current node onto the beginning of $w$, and then replace the current node by its parent.
\end{description}
These maps form a bijection (and its inverse) between $a$-valid sequences with at most $t$ digits and non-leftmost leaves below the $t^{th}$ supernode of $T$.  Furthermore, $v$ is the $\p{N+1}^{st}$ leaf from the left in $T$ if and only if the sequence is the $a$-Zeckendorf representation of $N$.
\end{prop}

Figure~\ref{fig:conbij} shows the bijection in Proposition~\ref{prop:ztree} for the sequence $a=\bk{2}$ (corresponding to the Conolly sequence).

\begin{figure}
\begin{center}
\resizebox{\textwidth-0.5in}{!}{
\tikzstyle{super}=[rectangle, draw]
\tikzstyle{regular}=[circle, draw]
\tikzstyle{leaf}=[ellipse, draw]
\begin{tikzpicture}[level distance = 2.2cm, grow=down]
\Tree[.{$\phantom{aaa}\iddots$}
[.\node[super]{$\phantom{0000}$};
  [.\node[super]{$\phantom{0000}$};
    [.\node[super]{$\phantom{0000}$};
      [.\node[super]{$\phantom{0000}$};
        [.\node[leaf]{$\phantom{0000}$};
        ]
        [.\node[leaf]{$\phantom{0000}$};
        ]
      ]
      \edge node[auto=left]{{\LARGE $1$}};
      [.\node[regular]{$\phantom{0000}$};
        [.\node[leaf]{$\phantom{0000}$};
        ]
        \edge node[auto=left]{{\LARGE $1$}};
        [.\node[leaf]{$\phantom{0}11\phantom{0}$};
        ]
      ]
    ]
    [.\node[regular]{$\phantom{0000}$};
      [.\node[regular]{$\phantom{0000}$};
        [.\node[leaf]{$\phantom{0000}$};
        ]
        [.\node[leaf]{$\phantom{0000}$};
        ]
      ]
      [.\node[regular]{$\phantom{0000}$};
        [.\node[leaf]{$\phantom{0000}$};
        ]
        [.\node[leaf]{$\phantom{0000}$};
        ]
      ]
    ]
  ]
  \edge node[auto=left]{{\LARGE $1$}};
  [.\node[regular]{$\phantom{0000}$};
   \edge node[auto=right]{{\LARGE $0$}};
    [.\node[regular]{$\phantom{0000}$};
      [.\node[regular]{$\phantom{0000}$};
        [.\node[leaf]{$\phantom{0000}$};
        ]
        [.\node[leaf]{$\phantom{0000}$};
        ]
      ]
      \edge node[auto=left]{{\LARGE $1$}};
      [.\node[regular]{$\phantom{0000}$};
       \edge node[auto=right]{{\LARGE $0$}};
        [.\node[leaf]{$1010$};
        ]
        [.\node[leaf]{$\phantom{0000}$};
        ]
      ]
    ]
    [.\node[regular]{$\phantom{0000}$};
      [.\node[regular]{$\phantom{0000}$};
        [.\node[leaf]{$\phantom{0000}$};
        ]
        [.\node[leaf]{$\phantom{0000}$};
        ]
      ]
      [.\node[regular]{$\phantom{0000}$};
        [.\node[leaf]{$\phantom{0000}$};
        ]
        [.\node[leaf]{$\phantom{0000}$};
        ]
      ]
    ]
  ]
]
\edge[draw=none]; {} ]
\end{tikzpicture}
}
\end{center}
\caption{An illustration of the bijection in Proposition~\ref{prop:ztree} for the digit sequences $11$ and $1010$ on the tree for the sequence $a=\bk{2}$.}
\label{fig:conbij}
\end{figure}
\begin{proof}
%We will construct a bijection between the set of valid words of length at most $i$ and the leaves below the $i^{th}$ supernode of $T$.  By Proposition~\ref{prop:zcount}, this will prove that there are $a_i-1$ leaves other than the leftmost leaf, or $a_i$ leaves in total.
%
%\begin{description}
%\item[Words to Leaves:] Let $w$ be a valid word of length $m\leq i$. Start at the $m^{th}$ supernode, and read $w$ from left to right.  As symbol $\ell$ is read from $w$, move from the current node to its \typel{\ell} child and make that the new current node.
%\item[Leaves to Words:] Let $v$ be a leaf below the $i^{th}$ supernode that is not the leftmost leaf.  Start with $v$ as the current node and $w$ as an empty word.  Until the current node is a supernode, concatenate the type of the current node onto the beginning of $w$, and then replace the current node by its parent.
%\end{description}

To show that the two operations define a bijection, we claim that they are both well-defined and are inverses of each other.
%, which will suffice to complete this proof.
\begin{description}
\item[Digits to Leaves is Well-Defined:] We must show that the final node we reach is a leaf and that it exists.  By Proposition~\ref{prop:tperf}, all the leaves of $T$ are at the same level, $M$ levels below the $M^{th}$ supernode.  So, if we can legally move to a child at each step, we end at a leaf. Proposition~\ref{prop:tchildcount} guarantees that all moves are possible.  We start at a supernode, so the first move is to one of its $\Lambda$ children (but not the leftmost one).  Thereafter, after having seen $j$ consecutive zeroes in the sequence, we are at a node $j$ levels below the most recent non-\typel{0} ancestor of the current node.  This node has $\Lambda_{j+1}$ children, which is the same number of distinct digits that are allowed after $j$ zeroes, by Proposition~\ref{prop:validz}.
\item[Leaves to Digits is Well-Defined:] We must show that $\dseq{d}{M}$ is an $a$-valid sequence with at most $t$ digits.  A supernode must be reached in at most $t$ steps, since the leaf we start at has the $t^{th}$ supernode as an ancestor.  So, $M\leq t$.  Since we do not start at the leftmost leaf and since we stop as soon as we reach a supernode, $d_M\neq0$.  Also, 
for each $i$, we have $0\leq d_i<\Lambda$, 
%$w\in\st{0,1,2,\ldots,\Lambda-1}^*$, 
since every node in $T$ has at most $\Lambda$ children.  Finally, we must verify that, if $d_i=\ell$ for some $i<M$, it is preceded by $\mu_\ell-1$ zeroes. Suppose $d_i=\ell$ for some $i<M$.  At the moment that digit was determined, the current node, call it $v_i$, was the root of a copy of $T_i$.  Since $i<M$, $v_i$'s parent was not a supernode. %, and there is a nonzero digit $d_{i'}$ with $i'>i$.  We can choose $i'$ to be the \emph{mimimal} index with this property. 
Hence, we can define $i'$ to be the minimal index such that $i'>i$ and $d_{i'}\neq0$. 
This means that $d_i$ is preceded by $i'-i-1$ zeroes.  At the moment that $d_{i'}$ was determined, the current node, call it $v_{i'}$, was the root of a copy of $T_{i'}$.  So, we have a copy of $T_{i'}$ that contains a copy of $T_i$ whose root is of \typel{\ell}.  This requires $\Lambda_{i'-i}\geq\ell+1$. % (see Figure~\ref{fig:mdiff}).  
This is equivalent to $\Lambda_{i'-i}>\ell$, or $\mu_\ell\leq i'-i$.  The right side of this last expression is one more than the number of zeroes preceding the digit $d_i$, so we have that the digit $\ell$ must be preceded by at least $\mu_\ell-1$ zeroes, as required.
\item[Inverses:] First, let $\dseq{d}{M}$ be an $a$-valid sequence with at most $t$ digits.  The first operation yields a leaf $v$.  
On the way to that leaf, each ancestor of $v$ from its most recent supernode ancestor to itself is visited, and the digits consist of the sequence of types of those nodes (except for the supernode).  Performing the second operation to pass from $v$ to an $a$-valid sequence $\dseq{d'}{M}$ constructs the sequence of types of $v$ and its non-supernode ancestors.  This is the same as $\dseq{d}{M}$, as required.

The argument that a non-leftmost leaf $v$ is preserved by passing to an $a$-valid sequence and back to a leaf is similar.
%Then, starting from $v$, the second operation yields a valid word $w'$.  Suppose for a contradiction that $w\neq w'$.  Let $m$ be the number of symbols from the end of the words where they last differ.  Suppose $w$ has a symbol $\ell_1$ in that position, and $w'$ has a symbol $\ell_2$ there.  By the first operation, $v$ must have its ancestor $m-1$ generations ago have Type~$\ell_1$.  But, then the $m^{th}$-to-last symbol in $w'$ would be $\ell_1$, which it 
\end{description}

We now claim that if $v'$ is a leaf to the right of $v$ (and both are below the $t^{th}$ supernode), then the sequence corresponding to $v'$ comes later in order than the sequence corresponding to $v$.  Let $\dseq{d'}{M}$ correspond to $v'$ and $\dseq{d}{M}$ to $v$.  Clearly we cannot have $M'<M$, as that would mean that $v'$ is a descendent of an earlier supernode than $v$ is (and hence it would be to the left of $v$).  So, suppose that $M=M'$, and consider the process of converting back from these sequences to leaves.  Let $i$ be the first index where $d_i\neq d_i'$.  At that moment, the current nodes for both conversions would be the same.  It must then be the case that the prime sequence moves to a farther-right node than the non-prime sequence does, as $v'$ is to the right of $v$.  In other words, $d_i'>d_i$, so the prime sequence is lexicographically later than the non-prime sequence, as required.

The claim that $v$ is the $\p{N+1}^{st}$ leaf from the left in $T$ if and only if $\dseq{d}{M}$ is the $a$-Zeckendorf representation of $N$ now follows by considering the facts that the bijection applies to every $a$-valid sequence and that a leaf farther to the right of another leaf corresponds to a sequence that comes later in order.
\end{proof}

%\begin{figure}
%\begin{center}
%%\includegraphics[width=300pt]{hof10000.eps}
%TODO
%\end{center}
%\caption{TODO}
%\label{fig:mdiff}
%\end{figure}

\begin{cor}\label{cor:ztree}
The number of leaves below the $t^{th}$ supernode of $T$ is $a_t$.
\end{cor}
\begin{proof}
Proposition~\ref{prop:ztree} gives a bijection between the set of $a$-valid sequences with at most $t$ digits and the non-leftmost leaves below the $t^{th}$ supernode of $T$.  By Proposition~\ref{prop:zcount}, this proves that there are $a_t-1$ leaves other than the leftmost leaf, or $a_t$ leaves in total.
\end{proof}

We now tie $a$-Zeckendorf representations back to leaf-counting functions.
\begin{prop}\label{prop:ztreecount}
Let $z_N$ denote the number of zeroes at the end of the $a$-Zeckendorf representation of the positive integer $N$, and let $s$ be a nonnegative integer. The number of times the value $N$ appears in the sequence $\seq{L_{T,s}\p{n}}_{n\geq1}$ is $1+z_N$, unless $N=a_i$ for some $i$, in which case it is $1+z_N+s$. 
\end{prop}
\begin{proof}
%(As usual, consider $s\geq0$.  We treat $s<0$ in Section~\ref{sec:neg}.)
%
%The value of $L_{T,s}\p{n}$ is defined as the number of leaves in $T\p{n,s}$.  (TODO define this notation)  
Let $v_N$ denote the $N^{th}$ leaf in $T$, and let $v_{N+1}$ denote the $\p{N+1}^{st}$ leaf.  Let $r$ denote the number of regular nodes between $v_N$ and $v_{N+1}$, and let $u$ denote the number of supernodes between them (always $0$ or $1$).  The number of times $N$ appears in $L_{T,s}$ is $r+1+su$.  Let $\dseq{d}{M}$ be the sequence obtained from $v_{N+1}$ under the bijection in Proposition~\ref{prop:ztree}.  The zeroes at the end of the sequence correspond precisely to the regular nodes in $T$ that are visited in preorder before $v_{N+1}$, so $r$ precisely equals $z_N$.  Furthermore, $\dseq{d}{M}$ is the $a$-Zeckendorf representation of $N$.  So, $N$ appears at least $1+z_N$ times.
%Let $n_t$ denote the label of the $t^{th}$ leaf, and let $n_{t+1}$ denote the label of the $\p{t+1}^{st}$ leaf.

The only situation in which $t$ could appear a different number of times than $1+z_N$ is if there is a supernode between $v_N$ and $v_{N+1}$.  This happens precisely when the $a$-Zeckendorf representations of $N-1$ and $N$ have different lengths.  By Corollary~\ref{cor:zcount}, this happens exactly when $N=a_i$ for some $i$.  In that case, there are $s$ extra labels, as required.
\end{proof}

\subsection{Leaf Counting Properties}
We now use the connection from Proposition~\ref{prop:ztree} between trees and $a$-Zeckendorf representations to explore additional properties of our sequences. To proceed, we study details of enumerating leaves in various subtrees of the trees $T_j$ and the skeleton $T$. A major outcome of this section is a pair of efficient algorithms for converting between $a$-Zeckendorf representations and numbers. Since the algorithms themselves are tangential to our discussion, we relegate them to Appendix~\ref{s:algs}.

%Note that there is a simple brute-force algorithm for going either direction: Starting from $1$, generate successive representations until generating the target number of representations (to go from numbers to representations) or until generating the target representation (to go from representations to numbers). This algorithm, while simple, is exponential in the length of the representation and linear in the number.

%Our algorithms make use of the 
%The bijection in Proposition~\ref{prop:ztree} is the main connection between trees and $a$-Zeckendorf representations. So, we first must explore leaf counts in various trees.
Recall (Definition~\ref{def:subtree}) that a subtree of a skeleton is a node along with all of its descendants. Similarly, for a nonnegative integer $j$, a subtree of a tree $T_j$ is a node along with all of its descendants. 
Let $\NN\p{j,t}$ denote the number of leaves in the subtree of $T_j$ rooted at the $t^{th}$ special node. 
We observe the following facts about the values of $\NN\p{j,t}$.%First, we have the following statements about leaf counts in subtrees of $T$.
\begin{lemma}\label{lem:Ntt}
For integers $j\geq0$,
\[
\NN\p{j,j}=\frac{a_{j+1}-a_j}{\Lambda-1}
\]
\end{lemma}
\begin{proof}
By Corollary~\ref{cor:ztree}, the subtree of $T$ rooted at its $\p{j+1}^{st}$ special node has $a_{j+1}$ leaves, and the subtree of $T$ rooted at its $j^{th}$ special node has $a_j$ leaves. So, there are $a_{j+1}-a_j$ leaves in $T$ that are below its $\p{j+1}^{st}$ special node and not below its $j^{th}$ special node. These leaves are partitioned equally into $\Lambda-1$ copies of $T_j$. Therefore,
\[
\NN\p{j,j}=\frac{a_{j+1}-a_j}{\Lambda-1},
\]
as required.
\end{proof}
\begin{lemma}\label{lem:Nrec}
For integers $0\leq t<j$, we have
\[
\NN\p{j,t}=\NN\p{j,t+1}-\pb{\Lambda_{j-t}-1}\pb{\frac{a_{t+1}-a_t}{\Lambda-1}}.
\]
\end{lemma}
\begin{proof}
Suppose $0\leq t<j$. Each leaf in the subtree of $T_j$ rooted at the $\p{t+1}^{st}$ special node is one of the following:
\begin{itemize}
\item in the subtree of $T_j$ rooted at the $t^{th}$ special node
\item in one of $\Lambda_{j-t}-1$ copies of $T_{t}$.
\end{itemize}
This gives rise to the recurrence
\[
\NN\p{j,t+1}=\NN\p{j,t}+\pb{\Lambda_{j-t}-1}\NN\p{t,t}=\NN\p{j,t}+\pb{\Lambda_{j-t}-1}\pb{\frac{a_{t+1}-a_{t}}{\Lambda-1}},
\]
with the last equality by Lemma~\ref{lem:Ntt}.
%\[
%\NN\p{t,t}=\frac{a_{t+1}-a_{t}}{\Lambda-1},
%\]
%which can be substituted to give
%\[
%\NN\p{j,t+1}=\NN\p{j,t}+\pb{\Lambda_{j-t}-1}\pb{\frac{a_{t+1}-a_{t}}{\Lambda-1}}.
%\]
Solving this expression for $\NN\p{j,t}$ yields the desired recurrence.
\end{proof}
\begin{prop}\label{prop:tjleafcount}
For integers $0\leq t\leq j$,
\[
\NN\p{j,t}=\frac{1}{\Lambda-1}\pb{a_{j+1}+\pb{\Lambda_{j-t}-1}a_t-\sum_{i=1}^{j-t}\lambda_ia_{j+1-i}}.
\]
\end{prop}
\begin{proof}
We proceed by induction on $j-t$. First, if $j-t=0$, we obtain
\[
\NN\p{j,j}=\frac{1}{\Lambda-1}\pb{a_{j+1}+\pb{0-1}a_j-0}=\frac{a_{j+1}-a_j}{\Lambda-1},
\]
which is the expression from Lemma~\ref{lem:Ntt}.

Now, suppose $j-t\geq1$ and that the formula holds when $t$ is replaced by $t+1$. Applying the recurrence from Lemma~\ref{lem:Nrec} yields
\begin{align*}
\NN\p{j,t}&=\NN\p{j,t+1}-\pb{\Lambda_{j-t}-1}\pb{\frac{a_{t+1}-a_t}{\Lambda-1}}\\
&=\frac{1}{\Lambda-1}\pb{a_{j+1}+\pb{\Lambda_{j-t-1}-1}a_{t+1}-\sum_{i=1}^{j-t-1}\lambda_ia_{j+1-i}}-\pb{\Lambda_{j-t}-1}\pb{\frac{a_{t+1}-a_t}{\Lambda-1}}\\
&=\frac{1}{\Lambda-1}\pb{a_{j+1}+\pb{\Lambda_{j-t-1}-1}a_{t+1}-\pb{\Lambda_{j-t}-1}a_{t+1}+\pb{\Lambda_{j-t}-1}a_{t}-\sum_{i=1}^{j-t-1}\lambda_ia_{j+1-i}}\\
&=\frac{1}{\Lambda-1}\pb{a_{j+1}-\lambda_{j-t}a_{t+1}+\pb{\Lambda_{j-t}-1}a_{t}-\sum_{i=1}^{j-t-1}\lambda_ia_{j+1-i}}\\
&=\frac{1}{\Lambda-1}\pb{a_{j+1}+\pb{\Lambda_{j-t}-1}a_{t}-\sum_{i=1}^{j-t}\lambda_ia_{j+1-i}},
\end{align*}
as required.
\end{proof}

Proposition~\ref{prop:tjleafcount} has two simple corollaries.
%\begin{cor}\label{cor:tjleafcount}
%The tree $T_j$ contains
%\[
%\frac{a_{j+1}-a_j}{\Lambda-1}
%\]
%leaves.
%\end{cor}
%\begin{proof}
%The tree $T_j$ itself is the subtree of $T_j$ rooted at the $j^{th}$ special node. The result follows from substituting $t=j$ into Item~\ref{it:tjleaves} of Proposition~\ref{prop:tleafcount}, as $\Lambda_0-1=-1$, and the summation is empty.
%\end{proof}

\begin{cor}\label{cor:tjfarleft}
If $j-t\geq k$, 
%the subtree of $T_j$ rooted at the $t^{th}$ special node has $a_t$ leaves.
$\NN\p{j,t}=a_t$.
\end{cor}
%\begin{proof}
%If $j-t\geq k$, the summation in 
%Proposition~\ref{prop:tjleafcount}
%%Item~\ref{it:tjleaves} of Proposition~\ref{prop:tleafcount} 
%equals $a_{j+1}$, by the recurrence $a$. The summation then cancels with the $a_{j+1}$ in the formula. Furthermore, $\Lambda_{j-t}=\Lambda$, so the coefficient on the $a_t$ cancels with the $\frac{1}{\Lambda-1}$ out front, leaving only the $a_t$.
%\end{proof}

\begin{cor}\label{cor:tj0}
%The subtree of $T_j$ rooted at the $0^{th}$ special node has $1$ leaf.
For all integers $j\geq0$, $\NN\p{j,0}=1$.
\end{cor}

Both Corollary~\ref{cor:tjfarleft} and Corollary~\ref{cor:tj0} are obvious. Corollary~\ref{cor:tjfarleft} also follows by observing that, if $j-t\geq k$, the subtree of $T_j$ rooted at its $t^{th}$ special node is isomorphic to the subtree of $T$ rooted at its $t^{th}$ special node. When deriving Corollary~\ref{cor:tjfarleft} from Proposition~\ref{prop:tjleafcount}, the key is that $\Lambda_{j-t}=\Lambda$ when $j-t\geq k$. Corollary~\ref{cor:tj0} is trivial, as the subtree in question is just a single leaf. When deriving Corollary~\ref{cor:tj0} from Proposition~\ref{prop:tjleafcount}, the expression inside the parentheses becomes $-1$ plus the sum of all $\lambda$ values. This value is $\Lambda-1$, which cancels with the $\frac{1}{\Lambda-1}$.

We now arrive at the following proposition, which extends Corollary~\ref{cor:tj0} and leads to a property of $a$-Zeckendorf representations that is analogous to a property of traditional Zeckendorf representations.
\begin{prop}\label{prop:Nsatrec}
For all integers $j\geq0$, $\NN\p{j,0}=1$, and for all integers $1\leq t\leq j$,
\[
\NN\p{j,t}=\sum_{i=1}^t\lambda_i\NN\p{j-i,t-i}+\sum_{i=t+1}^j\lambda_i.
\]
\end{prop}
\begin{proof}
The proof is by induction on $t$. If $t=0$, we have $\NN\p{j,0}$ by Corollary~\ref{cor:tj0}.  Now, let $t_0\geq0$, and suppose the formula holds for all values $t<t_0+1$.  By Lemma~\ref{lem:Ntt} and Lemma~\ref{lem:Nrec},
\[
\NN\p{j,t_0}=\NN\p{j,t_0+1}-\pb{\Lambda_{j-{t_0}}-1}\NN\p{t_0,t_0}.
\]
Rearranging gives
\begin{equation}\label{eq:Nrec}
\NN\p{j,t_0+1}=\NN\p{j,t_0}+\pb{\Lambda_{j-{t_0}}-1}\NN\p{t_0,t_0}.
\end{equation}
%Recall that
%\[
%\Lambda_{j-{t_0}}=\sum_{r=1}^{j-t_0}\lambda_r.
%\]
%This allows us to write
%\begin{equation}\label{eq:Nrec}
%\NN\p{j,t_0+1}=\NN\p{j,t_0}+\pb{-1+\sum_{r=1}^{j-t_0}\lambda_r}\NN\p{t_0,t_0}.
%\end{equation}
%Expanding the right side inductively gives
%\begin{align*}
%\NN\p{j,t_0+1}&=\sum_{i=1}^{t_0}\lambda_i\NN\p{j-i,t_0-i}+\sum_{i=t_0+1}^j\lambda_i+\pb{-1+\sum_{r=1}^{j-t_0}\lambda_r}\pb{\sum_{i=1}^{t_0}\lambda_i\NN\p{t_0-i,t_0-i}}\\
%&=
%\end{align*}
Now, consider the target expression
\[
\sum_{i=1}^{t_0+1}\lambda_i\NN\p{j-i,t_0+1-i}+\sum_{i=t_0+2}^j\lambda_i.
\]
We wish to use Equation~\ref{eq:Nrec} to rewrite this expression, but it only applies when the second argument is positive. So, we split the expression as
\[
\lambda_{t_0+1}\NN\p{j-i,0}+\sum_{i=1}^{t_0}\lambda_i\NN\p{j-i,t_0+1-i}+\sum_{i=t_0+2}^j\lambda_i.
\]
We now rewrite it as
\[
%\sum_{i=1}^{t_0+1}\lambda_i\pb{\NN\p{j-i,t_0-i}+\pb{-1+\sum_{r=1}^{j-t_0}\lambda_r}\NN\p{t_0-i,t_0-i}}+\sum_{i=t_0+2}^j\lambda_i.
\lambda_{t_0+1}+\sum_{i=1}^{t_0}\lambda_i\pb{\NN\p{j-i,t_0-i}+\pb{\Lambda_{j-t_0}-1}\NN\p{t_0-i,t_0-i}}+\sum_{i=t_0+2}^j\lambda_i.
\]
We now manipulate this expression to
%\begin{align*}
\[
%&\phantom{=}\sum_{i=1}^{t_0+1}\lambda_i\pb{\NN\p{j-i,t_0-i}+\pb{-1+\sum_{r=1}^{j-t_0}\lambda_r}\NN\p{t_0-i,t_0-i}}+\sum_{i=t_0+2}^j\lambda_i\\
%&=\sum_{i=1}^{t_0+1}\lambda_i\NN\p{j-i,t_0-i}-\sum_{i=1}^{t_0+1}\lambda_i\NN\p{t_0-i,t_0-i}+\sum_{i=1}^{t_0+1}\sum_{r=1}^{j-t_0}\lambda_i\lambda_r\NN\p{t_0-i,t_0-i}+\sum_{i=t_0+2}^j\lambda_i
%&\phantom{=}\sum_{i=1}^{t_0+1}\lambda_i\pb{\NN\p{j-i,t_0-i}+\pb{\Lambda_{j-t_0}-1}\NN\p{t_0-i,t_0-i}}+\sum_{i=t_0+2}^j\lambda_i\\
%&=
\lambda_{t_0+1}+\sum_{i=1}^{t_0}\lambda_i\NN\p{j-i,t_0-i}+\pb{\Lambda_{j-t_0}-1}\sum_{i=1}^{t_0}\lambda_i\NN\p{t_0-i,t_0-i}+\sum_{i=t_0+2}^j\lambda_i
%\end{align*}
\]
By induction, 
\[
\sum_{i=1}^{t_0}\lambda_i\NN\p{j-i,t_0-i}=\NN\p{j,t_0}-\sum_{i=t_0+1}^j\lambda_i
\]
and
\[
\sum_{i=1}^{t_0}\lambda_i\NN\p{t_0-i,t_0-i}=\NN\p{t_0,t_0}.
\]
So, our expression becomes
\[
\lambda_{t_0+1}+\NN\p{j,t_0}-\sum_{i=t_0+1}^j\lambda_i+\pb{\Lambda_{j-t_0}-1}\NN\p{t_0,t_0}+\sum_{i=t_0+2}^j\lambda_i.
\]
The initial $\lambda_{t_0+1}$ and the two remaining summations cancel out, leaving $\NN\p{j,t_0}+\pb{\Lambda_{j-t_0}-1}\NN\p{t_0,t_0}$, which equals $\NN\p{j,t_0+1}$ by Equation~\ref{eq:Nrec}, as required.
%the first summation here is\NN\p{j,t_0}-$, and the second summation is $\NN\p{t_0,t_0}$. So, we have
%\[
%\lambda_{t_0+1}+\NN\p{j,t_0}+\pb{\Lambda_{j-t_0}-1}\NN\p{t_0,t_0}+\sum_{i=t_0+2}^j\lambda_i.
%\]
%Combining the first term 
\end{proof}

Proposition~\ref{prop:Nsatrec} has the following useful corollary.
\begin{cor}\label{cor:Nsatrec}
For all integers $0\leq j\leq t$, the sequence $\seq{\NN\p{j+n,t+n}}_{n\geq0}$ satisfies the recurrence $a$, with the first $k$ terms as the initial condition.
\end{cor}
\begin{proof}
By Proposition~\ref{prop:Nsatrec},
\[
\NN\p{j+n,t+n}=\sum_{i=1}^{t+n}\lambda_i\NN\p{j+n-i,t+n-i}+\sum_{i=t+n+1}^{j+n}\lambda_i.
\]
If $n\geq k$, all terms in the second summation are necessarily zero, and all terms beyond $i=k$ in the first summation are zero as well. This yields
\[
\NN\p{j+n,t+n}=\sum_{i=1}^{k}\lambda_i\NN\p{j+n-i,t+n-i},
\]
which is precisely the recurrence $a$, as required.
\end{proof}

Corollary~\ref{cor:Nsatrec} leads to Proposition~\ref{prop:addzero}, which generalizes a key property of traditional Zeckendorf representations.

\begin{prop}\label{prop:addzero}
Let $a$ be a linear recurrence of order $k$, and let $\dseq{d}{M}$ be an $a$-valid sequence. Let $b_n$ be the number whose $a$-Zeckendorf representation is 
$\dseq{d^{(n)}}{M+n}$, where
%$d_0^{(n)},d_1^{(n)},\ldots,d_{M+n}^{(n)}$, where
\[
d_i^{(n)}=\begin{cases}
0 & i < n\\
d_{i-n} & i\geq n.
\end{cases}
\]
(The representation of $b_n$ is $\dseq{d}{M},0,0,\ldots,0$, with $n$ extra zeroes at the end.)
%$D_i=0,0,\ldots,0,d_0,d_1,\ldots,d_M$, where the number of digits in the representation $D_i$ is $i+M+1$ (so the number of extra zeroes is $i$). 
%$w0^i$ ($w$ followed by $i$ zeroes). 
The sequence $\seq{b_n}_{n\geq0}$ satisfies the recurrence $a$, with initial conditions $b_0,b_1,\ldots,b_{k-1}$.
\end{prop}

In order to prove Proposition~\ref{prop:addzero}, we need one more result. The following proposition allows us to prove Proposition~\ref{prop:addzero}, and it is the key to our algorithms in Appendix~\ref{s:algs}

%We conclude this subsection with one more proposition, which is the key to our algorithms in Appendix~\ref{s:algs} and is used in some other proofs later.
\begin{prop}\label{prop:Nform}
Let $\dseq{d}{N}$ be an $a$-valid sequence with $M+1$ digits that is the $a$-Zeckendorf representation of the number $N$. Let
\[
P=\st{(i,j):i<j\text{ and }d_i\neq0\text{ and }d_j\neq0\text{ and }d_p=0\text{ for }i<p<j}
\]
be the set of pairs of successive indices with nonzero digits. Then,
\[
N=a_M+\sum_{i:d_i\neq0}\pb{d_i-1}\NN\p{i,i}+\sum_{\p{i,j}\in P}\NN\p{j,i}.
\]
\end{prop}
\begin{proof}
We use the bijection from Proposition~\ref{prop:ztree} to convert from $a$-valid sequences to tree paths. By Proposition~\ref{prop:ztree} the representation $\dseq{d}{M}$ describes a traversal down the skeleton $T$ for recurrence $a$ to a leaf. Then, $N$ is the number of non-leftmost leaves at or to the left of that leaf. Equivalently, $N$ is the number of leaves to the left of that leaf. First, recall that $d_M\neq0$. So, the first step in the path is from the $M^{th}$ supernode to a regular node or leaf. To the left of this step is the left subtree of $T$ rooted at the $\p{M-1}^{st}$ supernode, which has $a_M$ leaves, by Corollary~\ref{cor:ztree}. 

At each subsequent step, processing a digit $d_i$ with $i<M$, one of the following cases happens:
\begin{description}
\item[$d_i=0$:] In this, case there are no additional leaves to the left of the current node than there were at the prior step. So, a zero introduces no terms into the formula for $N$.
\item[$d_i=\ell>0$:] This case involves a step from a node $v_i$ to its \typel{\ell} child, and it does introduce additional leaves. First, the leaves of the left subtree of $v_i$ are now to the left of the current node. These are all of the leaves below the $i^{th}$ special node of $T_j$, where $j$ is the level where we most recently moved to a non-\typel{0} node. The number of such leaves is $\NN\p{j,i}$, and $\p{i,j}\in P$. Second, the leaves below all \typel{\ell'} children of $v$ are now to the left of the current node, for each $0<\ell'<\ell$. There are $\ell-1$ such children of $v_i$, and each is the root of a copy of $T_i$. Hence, each of these trees has $\NN\p{i,i}$ leaves.
\end{description}
Starting from $a_M$ and adding up the contributions of all of the nonzero digits gives the required formula for $N$.
%Algorithm~\ref{alg:zn} works by starting from $0$ and accumulating the value of $N$ as $w$ is traversed from left to right. As we traverse $w$, we travel down $T$ until we reach a leaf at the end. There are two ways that leaves can appear to the left of the leaf reached along this path:
%On Line~\ref{ln:znquadam} of Algorithm~\ref{alg:znquad}, we initialize $N$ to $a_m$ to account for these leaves from the start. Aside from these $a_m$ leaves, there are two other ways leaves can appear to the left of the leaf reached along the path:
%\begin{itemize}
%\item On a step from a node $v$ to a non-\typel{0} child, all leaves in the subtree rooted at $v$'s \typel{0} child are to the left of the reached node.
%\item On a step from a node $v$ to a \typel{\ell} child where $\ell>0$, all leaves in the subtrees rooted at children of $v$ of type greater than $0$ and less than $\ell$ are to the left of the reached node.
%\end{itemize}
%Any step to a non-\typel{0} child is to the root of some $T_j$. Line~\ref{ln:znquadleft} accumulates the leaf counts of subtrees rooted at \typel{0} nodes, and Line~\ref{ln:znquadnotleft} accumulates the other leaf counts. Algorithm~\ref{alg:znquad} is written traversing the path from bottom up, though it could have also been implemented top-down.
%Aside from the initial accumulation of $a_m$, each subsequent accumulation of the first type involves adding a value $\NN\p{j,t}$ for some $0\leq t<j$. Accumulations of the second type involve adding values of the form $\NN\p{j,j}$.
\end{proof}

We now conclude this subsection by proving Proposition~\ref{prop:addzero}.
\begin{proof}
%Suppose $w$ has length $m+1$.
By Proposition~\ref{prop:Nform}, $b_0$ can be written as a linear combination of $a_M$ and values $\NN\p{j,t}$ where one of the following holds:
\begin{itemize}
\item $j=t$ and $d_t\neq0$
\item $j>t$, $d_j\neq0$, $d_t\neq0$, and all digits between $d_j$ and $d_t$ are zero.
\end{itemize}
Note that, for each $n$, $d_t\neq0$ if and only if $d_{t+n}^{(n)}\neq0$. This means that $b_n$ can be written as a linear combination of $a_{M+n}$ and values $\NN\p{j,t}$ where one of the following holds:
\begin{itemize}
\item $j=t$ and $d_{t-n}\neq0$
\item $j>t$, $d_{j-n}\neq0$, $d_{t-n}\neq0$, and all digits between $d_{j-n}$ and $d_{t-n}$ are zero.
\end{itemize}
The key observation is that the coefficient on $a_M$ for $b_0$ is the same as the coefficient on $a_{M+n}$ for $b_n$, and any coefficient on $\NN\p{j,t}$ for $b_0$ is the same as the coefficient on $\NN\p{j+n,t+n}$ for $b_n$. By definition, the $a$ sequence satisfies the $a$ recurrence. By Corollary~\ref{cor:Nsatrec}, so do the numbers $\seq{\NN\p{j+n,t+n}}_{n\geq0}$. Therefore, the numbers $\seq{b_n}_{n\geq0}$ satisfy the $a$ recurrence, as required.
\end{proof}

%\subsection{Other Properties}\label{ss:op}

\section{Asymptotics}\label{sec:asympt}

In this subsection, we study the asymptotic behavior of the sequences $L_{T,s}$.  
Our analysis depends on a result known sometimes as Ostrovsky's Theorem or sometimes as the Cauchy-Ostrovsky Theorem.
\begin{theorem}\cite{prasolovpolys}\label{thm:ostro}
Let $t$ be a positive integer, let $\alpha_1,\alpha_2,\ldots,\alpha_t$ be nonnegative real numbers with at least one of them positive.  Let
\[
p\p{x}=x^t-\sum_{i=1}^t\alpha_ix^{t-i}.
\]
If the greatest common divisor of the indices $i$ such that $\alpha_i>0$ is $1$, then $p\p{x}$ has a unique positive real root $\kappa$ (counted with multiplicity), and that root is the unique root of $p\p{x}$ of largest modulus.
\end{theorem}

We use the following lemma:
\begin{lemma}\label{lem:poly}
Let the sequence $a=\bk{\lambda_1,\lambda_2,\ldots,\lambda_k}$ be as described in Subsection~\ref{ss:not}. %(TODO require not k=1, l=1).  
%Let $\kappa$ be the positive real root of the polynomial
Consider the polynomial
\[
p\p{x}=x^k-\sum_{i=1}^k\lambda_ix^{k-i}.
\]
%Let $\kappa=\kappa_1,\kappa_2,\kappa_3,\ldots,\kappa_k$ be the (complex) roots of $p\p{x}$ in decreasing order of modulus (listed with multiplicity, ties broken arbitrarily). 
Let $\kappa$ be a (complex) root of $p\p{x}$ of maximum modulus. The number $\kappa$ is unique, not a multiple root, real, greater than $1$, and the only positive real root of $p\p{x}$.
%that lies in $\p{1,\infty}$ (as per Lemma~\ref{lem:poly}).  %If the modulus of each other root of $p\p{x}$ is at most $1$, then for 
%The number $\kappa$ is the unique positive real root of $p\p{x}$. 
%Also, $\kappa>1$, and $\ab{\kappa}>\ab{\kappa_2}$. 
%Furthermore, %the sequence $a$ satisfies
%%\[
%%\limi{n}{\frac{a_n}{a_{n-1}}}=\kappa
%%\]
%there are nonnegative integers $e_1,e_2,\ldots,e_k$ with $e_1$ necessarily equal to $0$ and (complex) constants $B_1,B_2,\ldots,B_k$ with $B_1$ real and positive such that
%\[
%a_n=\sum_{i=1}^kB_in^{e_i}\kappa_i^n.
%\]
Furthermore, there exists a positive constant $B$ such that $a_n=B\kappa^n+o\p{\kappa^n}$.
\end{lemma}
\begin{proof}
We always have $\lambda_1>0$ and $\lambda_k>0$. Since $\gcd\p{1,k}=1$, by Theorem~\ref{thm:ostro} $p\p{x}$ has a unique positive real root, and it the unique root of largest modulus among all roots of $p\p{x}$. This root is therefore $\kappa$. We now claim that $\kappa>1$. Note that $p\p{1}=1-\Lambda$. Since $\Lambda\geq2$, we have $p\p{1}<0$. Also, because $p\p{x}$ has positive leading coefficient,
\[
\limi{x}{p\p{x}}=\infty.
\]
So, by the Intermediate Value Theorem, $p\p{x}$ has a root in $\p{1,\infty}$, which must be $\kappa$, as required.

Now, we must prove that there is a positive constant $B$ such that $a_n=B\kappa^n+o\p{\kappa^n}$. The polynomial $p\p{x}$ is the characteristic polynomial of the recurrence $a$, so classical theory of linear recurrences
 and the fact that $\kappa$ is not a multiple root, 
tells us that 
%there are nonnegative integers $e_1,e_2,\ldots,e_k$ and (complex) constants $B_1,B_2,\ldots,B_k$ such that
%\[
%a_n=\sum_{i=1}^kB_in^{e_i}\kappa_i^n.
%\]
there exists a constant $B$ such that $a_n=B\kappa^n+o\p{\kappa^n}$. 
So, we need only argue that 
%$e_1=0$ and $B_1>0$. 
$B>0$. %Since $\ab{\kappa}>\ab{\kappa_2}$, the root $\kappa$ has multiplicity $1$. An exponent $e_i$ can only be positive if $\kappa_i$ is a multiple root, so $e_1=0$ necessarily.
First, we claim that $B\neq0$. Recall that we have $a_{1-k}=a_{2-k}=\cdots=a_0=1$, and $a_1=\Lambda>1$. Suppose for a contradiction that $B=0$. Consider the recurrence relation $b$ of order $k-1$ with characteristic polynomial $\frac{p\p{x}}{x-\kappa}$. This recurrence has the same characteristic roots as $a$, except it does not have $\kappa$ as a root. Since $B=0$, our sequence $\seq{a_n}$ also satisfies the recurrence $b$. Since $b$ has order $k-1$, $k-1$ initial conditions are sufficient to determine the sequence. The first $k-1$ terms are $a_{1-k}=a_{2-k}=\cdots=a_{-1}=1$. We know that $a_0=1$ also. Now, to compute $a_1$, we need only the preceding $k-1$ terms, which are all $1$. Since the same values were used to compute $a_0$, we must have $a_1=1$, a contradiction. Therefore, $B\neq0$. It immediately follows that $B>0$, as if $B$ were negative the sequence would eventually contain negative terms.

%Now, it must be the case that $B_1\in\mathbb{R}$, as otherwise we would not obtain real numbers in our sequence. So, we need only show that $B_1>0$. First, we claim that $B_1\neq0$. Recall that we have $a_{1-k}=a_{2-k}=\cdots=a_0=1$, and $a_1=\Lambda>1$. Suppose for a contradiction that $B_1=0$. Consider the recurrence relation $b$ of order $k-1$ with characteristic polynomial $\frac{p\p{x}}{x-\kappa}$. This recurrence has the same characteristic roots as $a$, except it does not have $\kappa$ as a root. Since $B_1=0$, our sequence $\seq{a_n}$ also satisfies the recurrence $b$. Since $b$ has order $k-1$, $k-1$ initial conditions are sufficient to determine the sequence. The first $k-1$ terms are $a_{1-k}=a_{2-k}=\cdots=a_{-1}=1$. We know that $a_0=1$ also. Now, to compute $a_1$, we need only the preceding $k-1$ terms, which are all $1$. Since the same values were used to compute $a_0$, we must have $a_1=1$, a contradiction. Therefore, $B_1\neq0$. It immediately follows that $B_1>0$, as otherwise the sequence would eventually contain negative terms.
\end{proof}
%TODO maybe only care about $a_n=B\kappa^n+o\p{\kappa^n}$

We also need the following:
%\begin{defin}\label{def:capa}
%Define the sequence $\seq{A_n}_{n\geq0}$ as the sequence of partial sums of $\seq{a_n}_{n\geq0}$, i.e.
%\[
%A_n=\sum_{i=0}^na_i.
%\]
%\end{defin}
\begin{lemma}\label{lem:Aasym}
Let $\kappa>1$, and let $\seq{b_n}_{n\geq0}$ be a sequence such that $b_n=B\kappa^n+o\p{\kappa^n}$ for some $B>0$. Define the sequence $\seq{B_n}_{n\geq0}$ as the sequence of partial sums of $\seq{b_n}_{n\geq0}$, i.e.
\[
B_n=\sum_{i=0}^nb_i.
\]
%Let the sequence $a=\bk{\lambda_1,\lambda_2,\ldots,\lambda_k}$ be as described in Subsection~\ref{ss:not}, and let $\kappa$ be the number from Lemma~\ref{lem:poly} such that $a_n=B\kappa^n+o\p{\kappa^n}$ for some positive constant $B$. 
We have $B_n=B\pb{\frac{\kappa}{\kappa-1}}\kappa^n+o\p{\kappa^n}$.
\end{lemma}
\begin{proof}
We have, using Lemma~\ref{lem:poly},
\begin{align*}
B_n&=\sum_{i=0}^nb_i\\
&=\sum_{i=0}^n\pb{B\kappa^i+o\p{\kappa^i}}\\
&=B\sum_{i=0}^n\kappa^i+o\!\pb{\sum_{i=0}^n\kappa^i}\\
&=B\pb{\frac{\kappa^{n+1}-1}{\kappa-1}}+o\!\pb{\frac{\kappa^{n+1}-1}{\kappa-1}}\\
&=B\pb{\frac{\kappa^{n+1}}{\kappa-1}}+o\p{\kappa^n}\\
&=B\pb{\frac{\kappa}{\kappa-1}}\kappa^n+o\p{\kappa^n},
\end{align*}
as required.
\end{proof}
%\begin{lemma}\label{lem:lincomblim}
%Let the sequence $a=\bk{\lambda_1,\lambda_2,\ldots,\lambda_k}$ be as described in Subsection~\ref{ss:not}, and let $\kappa$ be the number from Lemma~\ref{lem:poly} such that $a_n=B\kappa^n+o\p{\kappa^n}$ for some positive constant $B$. Let $\seq{\nu_n}_{n\geq0}$ be a sequence of nonnegative integers where at least one is positive, and let $c$ be a nonnegative integer.  Then,
%\[
%\limi{n}{\frac{\displaystyle\sum_{i=0}^n\nu_ia_i}{c+\displaystyle\sum_{i=0}^n\nu_iA_i}}=1-\frac{1}{\kappa}.
%\]
%\end{lemma}
%\begin{proof}
%Note that both the numerator and denominator are positive, because the coefficients are We have (using Lemma~\ref{lem:Aasym}),
%\begin{align*}
%\frac{\displaystyle\sum_{i=0}^n\nu_ia_i}{c+\displaystyle\sum_{i=0}^n\nu_iA_i}&=\frac{\displaystyle\sum_{i=0}^n\nu_i\pb{B\kappa^i+o\p{\kappa^i}}}{c+\displaystyle\sum_{i=0}^n\nu_i\pb{B\pb{\frac{\kappa}{\kappa-1}}\kappa^i+o\p{\kappa^i}}}\\
%&=\frac{B\displaystyle\sum_{i=0}^n\nu_i\kappa^i+o\pb{\sum_{i=0}^n\nu_i\kappa^i}}{c+B\displaystyle\pb{\frac{\kappa}{\kappa-1}}\displaystyle\sum_{i=0}^n\nu_i\kappa^i+o\pb{\sum_{i=0}^n\nu_i\kappa^i}}.
%\end{align*}
%Now, we divide the numerator and denominator by
%\[
%TODO
%\]
%\end{proof}
\begin{lemma}\label{lem:leafasym}
Let $a$ and $\kappa$ be as in Lemma~\ref{lem:poly}, and let $B$ be as in Lemma~\ref{lem:Aasym}. Let $t$ be a positive integer, and fix a nonnegative integer $c$. There exists a constant $B_{t,c}$ such that $L\p{t+c,t}=B_{t,c}\cdot\kappa^t+o\p{\kappa^t}$.
%For all integers $0<t<j$, there exists a constant $B_{j,t}$ such that $\NN\p{j,t}=B_{j,t}\cdot\kappa^t+o\p{\kappa^t}$, where the quantity $j-t$ is treated as fixed.
\end{lemma}
\begin{proof}
There are three cases to consider:
\begin{description}
\item[$c=0$:] In this case, $\NN\p{t+c,t}=\NN\p{t,t}=\frac{a_{t+1}-a_t}{\Lambda-1}$ by Lemma~\ref{lem:Ntt}. We have
\begin{align*}
\NN\p{t,t}&=\frac{1}{\Lambda-1}\pb{B\cdot\kappa^{t+1}-B\cdot\kappa^t}+o\p{\kappa^t}\\
&=\frac{B\pb{\kappa-1}}{\Lambda-1}\kappa^t+o\p{\kappa^t}.
\end{align*}
Hence, taking $B_{t,0}=\frac{B\pb{\kappa-1}}{\Lambda-1}$ suffices.
\item[$c\geq k$:] In this case, $\NN\p{t+c,t}=a_t$ by Corollary~\ref{cor:tjfarleft}. So, taking $B_{t+c,t}=B$ suffices to make $\NN\p{t+c,t}=B\cdot\kappa^t+o\p{\kappa^t}$.
\item[$1\leq c<k$:] By Proposition~\ref{prop:tjleafcount},
\[
\NN\p{t+c,t}=\frac{1}{\Lambda-1}\pb{a_{t+c+1}+\pb{\Lambda_{c}-1}a_t-\sum_{i=1}^{c}\lambda_ia_{t+c+1-i}}.
\]
We know that
\[
a_{t+c+1}=\sum_{i=1}^k\lambda_ia_{t+c+1-i},
\]
so we can rewrite the above expression as
\[
\NN\p{t+c,t}=\frac{1}{\Lambda-1}\pb{\pb{\Lambda_{c}-1}a_t+\sum_{i=c+1}^k\lambda_ia_{t+c+1-i}}.
\]
We now manipulate the expression:

\begin{align*}
\NN\p{t+c,t}&=\frac{1}{\Lambda-1}\pb{\pb{\Lambda_{c}-1}a_t+\sum_{i=c+1}^k\lambda_ia_{t+c+1-i}}\\
&=\frac{1}{\Lambda-1}\pb{\pb{\Lambda_{c}-1}B\cdot\kappa^t+\sum_{i=c+1}^k\lambda_iB\cdot\kappa^{t+c+1-i}}+o\p{\kappa^t}\\
&=\frac{1}{\Lambda-1}\pb{\pb{\Lambda_{c}-1}B+\sum_{i=c+1}^k\lambda_iB\cdot\kappa^{c+1-i}}\kappa^t+o\p{\kappa^t}.
\end{align*}
This is in the required form; we can take
\[
B_{t,c}=\frac{1}{\Lambda-1}\pb{\pb{\Lambda_{c}-1}B+\sum_{i=c+1}^k\lambda_iB\cdot\kappa^{c+1-i}}.
\]
\end{description}
%\begin{description}
%\item[$t=j$:] In this case, $\NN\p{j,t}=\frac{a_{j+1}-a_j}{\Lambda-1}$ by Lemma~\ref{lem:Ntt}. We have
%\begin{align*}
%\NN\p{j,t}&=\frac{1}{\Lambda-1}\pb{B\cdot\kappa^{j+1}-B\cdot\kappa^j}+o\p{\kappa^j}\\
%&=\frac{B\pb{\kappa-1}}{\Lambda-1}\kappa^j+o\p{\kappa^j}.
%\end{align*}
%Since $j=t$, taking $B_{j,j}=\frac{B\pb{\kappa-1}}{\Lambda-1}$ suffices.
%\item[$j-t\geq k$:] In this case, $\NN\p{j,t}=a_t$ by Corollary~\ref{cor:tjfarleft}. So, taking $B_{j,t}=B$ suffices to make $\NN\p{j,t}=B\cdot\kappa^t+o\p{\kappa^t}$.
%\item[$1\leq j-t<k$:] By Proposition~\ref{prop:tjleafcount},
%\[
%\NN\p{j,t}=\frac{1}{\Lambda-1}\pb{a_{j+1}+\pb{\Lambda_{j-t}-1}a_t-\sum_{i=1}^{j-t}\lambda_ia_{j+1-i}}.
%\]
%We know that
%\[
%a_{j+1}=\sum_{i=1}^k\lambda_ia_{j+1-i},
%\]
%so we can rewrite the above expression as
%\[
%\NN\p{j,t}=\frac{1}{\Lambda-1}\pb{\pb{\Lambda_{j-t}-1}a_t+\sum_{i=j-t+1}^k\lambda_ia_{j+1-i}}.
%\]
%We now manipulate the expression:
%
%\begin{align*}
%\NN\p{j,t}&=\frac{1}{\Lambda-1}\pb{\pb{\Lambda_{j-t}-1}a_t+\sum_{i=j-t+1}^k\lambda_ia_{j+1-i}}\\
%&=\frac{1}{\Lambda-1}\pb{\pb{\Lambda_{j-t}-1}B\cdot\kappa^t+\sum_{i=j-t+1}^k\lambda_iB\cdot\kappa^{j+1-i}}+o\p{\kappa^t}\\
%&=\frac{1}{\Lambda-1}\pb{\pb{\Lambda_{j-t}-1}B\cdot+\sum_{i=j-t+1}^k\lambda_iB\cdot\kappa^{j+1-i-t}}\kappa^t+o\p{\kappa^t},
%\end{align*}
%which is of the required form.
%\end{description}
\end{proof}

Our main result on asymptotics is the following:
\begin{theorem}\label{thm:asympt}
Let the sequence $a=\bk{\lambda_1,\lambda_2,\ldots,\lambda_k}$ be as described in Subsection~\ref{ss:not}. %(TODO require not k=1, l=1).  
%Let $\kappa$ be the positive real root of the polynomial
Consider the polynomial
\[
p\p{x}=x^k-\sum_{i=1}^k\lambda_ix^{k-i}.
\]
%that lies in $\p{1,\infty}$ (as per Lemma~\ref{lem:poly}).  %If the modulus of each other root of $p\p{x}$ is at most $1$, then for 
This polynomial has a unique real root $\kappa>1$, which is its unique root of largest modulus, and for 
%For 
any nonnegative integer $s$, the sequence $\seq{L_{T,s}\p{n}}_{n\geq1}$ has the property that
\[
\limi{n}{\frac{L_{T,s}\p{n}}{n}}=1-\frac{1}{\kappa}.
\]
%Otherwise, this limit does not exist.
\end{theorem}

%We use the following lemma:

%We can now prove Theorem~\ref{thm:asympt}.
\begin{proof}

Let $N$ be a positive integer, and let $n_N$ denote the smallest value such that $L_{T,s}\p{n_N}=N$. The number $N$ has an $a$-Zeckendorf representation, and, by Proposition~\ref{prop:Nform},
\[
N=a_M+\sum_{i:d_i\neq0}\pb{d_i-1}\NN\p{i,i}+\sum_{\p{i,j}\in P}\NN\p{j,i},
\]
where $P$ is as in the statement of Proposition~\ref{prop:Nform}. Note that $a_M$ is the number of leaves in the tree rooted at the $M^{th}$ supernode of the skeleton $T$ defined by $a$. So, every term in this sum is the number of leaves in some tree rooted at some special node of some tree $T$ or $T_j$. This formula puts these leaves together to count the first $N$ leaves of $T$. The label of this $N^{th}$ leaf is then $n_N$.%Traversing the path in $T$ determined by the $a$-Zeckendorf representation of $N$ and by Proposition~\ref{prop:ztree} leads to the $\p{N+1}^{st}$ leaf in $T$. The label of the $N^{th}$ leaf is then $n_N$.%$p$ of this leaf is the smallest index where $L_{T,s}=N+1$, so $a_N=p-1$.

Let $N^*$ denote the number of \emph{nodes} preceding or at the $N^{th}$ leaf of $T$. Recall that deleting the leaves from $T_j$ results in $T_{j-1}$ and deleting the leaves from $T$ results in $T$ (Propositions~\ref{prop:tilr} and~\ref{prop:tlr}). So, the number of nodes in the  tree rooted at the $M^{th}$ supernode of $T$ is
\[
\sum_{r=0}^Ma_r,
\]
and the number of nodes in the subtree rooted at the $t^{th}$ special node of $T_j$ is
\[
\sum_{r=0}^i\NN\p{j-r,i-r}.
\]
Putting this together with the formula for $N$ yields
\[
N^*=\sum_{r=0}^Ma_r+\sum_{i:d_i\neq0}\pb{d_i-1}\sum_{r=0}^i\NN\p{i-r,i-r}+\sum_{\p{i,j}\in P}\sum_{r=0}^i\NN\p{j-r,i-r}.
\]

Let $z$ be the number of zeroes at the end of the $a$-Zeckendorf representation of $N$. By Proposition~\ref{prop:ztreecount}, the largest index $t$ such that $L_{T,s}\p{t}$ equals $N$ is $n_N+z$, unless $N=a_M$, in which case it is $n_N+z+s=n_N+m+s$. Note that it must be the case that $z\leq M$.

Let us now examine the asymptotics of $N$ and $N^*$.  By Lemmas~\ref{lem:poly} and~\ref{lem:leafasym},
\[
N=B\cdot\kappa^M+\sum_{i:d_i\neq0}\pb{d_i-1}B_{i,0}\cdot\kappa^i+\sum_{\p{i,j}\in P}B_{i,j-i}\cdot\kappa^i+o\p{\kappa^M},
\]
where each $B$ term is a constant. Combining the aforementioned lemmas with Lemma~\ref{lem:Aasym} yields
\[
N^*=B\pb{\frac{\kappa}{\kappa-1}}\cdot\kappa^M+\sum_{i:d_i\neq0}\pb{d_i-1}B_{i,0}\pb{\frac{\kappa}{\kappa-1}}\cdot\kappa^i+\sum_{\p{i,j}\in P}B_{i,j-i}\pb{\frac{\kappa}{\kappa-1}}\cdot\kappa^i+o\p{\kappa^M}.
\]
In other words, $N^*=\pb{\frac{\kappa}{\kappa-1}}N+o\p{\kappa^M}$.

Let us now consider ratios $\frac{N}{n}$, where $L_{T,s}\p{n}=N$. We know that such an $n$ must satisfy $n_N\leq n\leq n_N+M+s=n_N+O\p{M}$
%Since $m+1$ is the height of the tree used to find the $N^{th}$ leaf, we have $m=O\p{\log N}$.

The number $N^*$ is closely related to number $n_N$, and they are, in fact, equal, if each node receives exactly one label. But, there are $M+1$ supernodes preceding the $N^{th}$ leaf ($M$ supernodes if $N=a_M$), each of which has $s$ labels. So, if $N\neq a_M$, we have $n_N=N^*+\pb{s-1}\pb{M+1}$, and if $N=a_M$, we have $n_N=N^*+\pb{s-1}M$. In any case, 
%$n_N=N^*+O\p{\log N}$.
$n_N=N^*\pm O\p{M}$. 
Going forward, note that $M=o\p{\kappa^M}$.

Let us now examine the extreme asymptotics, i.e. the upper bound of $\frac{N}{N^*-O\p{M}}$ and the lower bound of $\frac{N}{N^*+O\p{M}}$. We first have
\begin{align*}
\frac{N}{n}&\leq\frac{N}{N^*-O\p{M}}\\
&=\frac{N}{\pb{\frac{\kappa}{\kappa-1}}N+o\p{\kappa^M}-O\p{M}}\\
&=\frac{1}{\pb{\frac{\kappa}{\kappa-1}}+o\!\pb{\frac{\kappa^M}{N}}}\\
&=\frac{\kappa-1}{\kappa+o\p{1}}.
%\\
%&=1-\frac{1}{\kappa}+o\p{1}.
\end{align*}
We also have
\[
\frac{N}{n}\geq\frac{N}{N^*+O\p{M}},
\]
but the calculation is essentially the same, since $\kappa^M$ dominates $M$. So, we have
\[
\lim_{n\to\infty}\frac{L_{T,s}\p{n}}{n}=1-\frac{1}{\kappa},
\]
as required.

\end{proof}

%\subsection{Generalized Version of Theorem~\ref{thm:main}}

%\subsection{Limitations}

\section{Future Work}\label{sec:fut}
%lambda_1=0, negative coefficients, more generalizations with negative coefs, negative labels in non-supernodes, other trees
This work describes, for each homogeneous linear recurrence relation
\[
a_n=\sum_{i=1}^k\lambda_ia_{n-i}
\]
with nonnegative coefficients and $\lambda_1>0$, an infinite family of slow sequences satisfying nested recurrence relations based on the recurrence $a$. One question that remains is whether a similar result can be obtained when $\lambda_1=0$. The construction of a nested recurrence in Theorem~\ref{thm:main} still works for these for such recurrences, though the corresponding tree construction fails. Experimentally, it appears that slow solutions do exist for the nested recurrences constructed from relations with $\lambda_1=0$, but their connection to trees is unclear. For example, the recurrence $a=\bk{0,2}$ results in the nested recurrence
\[
C\p{n}=C(n-C(n-1)-C(n-1-C(n-1)))+C(n-1-C(n-2)-C(n-2-C(n-2))).
\]
The initial conditions $C\p{1}=1$, $C\p{2}=1$, $C\p{3}=2$ appear to result in a slow solution where each term appears the same number of times it appears in the Conolly sequence, except that the number $2^i$ appears $2^i+i+1$ times. Interestingly, this is the leaf counting function of the Conolly skeleton if the $i^{th}$ supernode receives $2^i$ labels. Similarly, the recurrence $a=\bk{0,1,1}$ results in the nested recurrence
\begin{align*}
&C(n-C(n-1)-C(n-C(n-1)-1))\\
&\phantom{C}+C(n-1-C(n-2)-C(n-2-C(n-2))\\
&\phantom{C+}-C(n-2-C(n-2)-C(n-2-C(n-2)))).
\end{align*}
The initial conditions $C\p{1}=1$, $C\p{2}=1$, $C\p{3}=2$, $C\p{4}=2$, $C\p{5}=2$ appear to result in a slow solution.  This solution appears to be connected to the sequence generated by the recurrence $a$ ($1,2,2,3,4,5,7,9,12,16,\ldots$, a shift of the Padovan sequence A000931), but the precise connection is not obvious. In general, members of this sequence seem to appear with higher frequencies than other nearby terms, with all observed record frequencies coming from terms in this sequence.

It may also be possible to obtain nested recurrences with slow solutions related to linear recurrence relations with some negative coefficients. With negative coefficients, the construction in Theorem~\ref{thm:main} no longer works, but perhaps it can be adapted in some way.

Finally, there are other sequences besides those related to Conolly's sequence that satisfy nested recurrences and are connected to trees~\cites{isgurthesis,nonslow,nonmonotgol,golomb}, such as Golomb's sequence $G\p{n}=G\p{n-G\p{n-1}}+1$, $G\p{1}=1$. It is likely possible to generalize many of these constructions to use multiple levels of pruning to obtain a generalization of the results in this paper.

%TODO NEGATIVE SHIFTS ARE IN~\cite{isgur2}!!!

\bibliography{bibliography.bib}

\appendix

\section{Some Examples}\label{sec:examples}

In this appendix, we examine a few examples that summarize all of our major results.
\subsection{Recurrences $a_n=\lambda\cdot a_{n-1}$}
We have already seen the example of the Conolly sequence, which arises from the recurrence $a_n=2a_{n-1}$ and results in a perfect binary tree. More generally, taking $\lambda\geq2$ results in the following:
\begin{itemize}
\item The tree $T_j$ is a perfect $\lambda$-ary tree with $j+1$ levels. (Definition~\ref{def:ti})
\item The skeleton $T$ is a perfect $\lambda$-ary tree. (Definition~\ref{def:t})
\item The $a$-Zeckendorf representation of $N$ is the base-$\lambda$ representation of $N$. (Proposition~\ref{prop:placeval})
\item The number of occurrences of the number $N$ in the sequence $\seq{L_{T,s}\p{n}}_{n\geq1}$ is one more than the number of zeroes at the end of the base-$\lambda$ representation of $N$, unless $N=\lambda^t$ for some $t$, in which case the number of occurrences is $t+1+s$. (Proposition~\ref{prop:ztreecount})
\item The sequence $\seq{L_{T,s}\p{n}}_{n\geq1}$ has the property that
\[
\lim_{n\to\infty}\frac{L_{T,s}\p{n}}{n}=1-\frac{1}{\lambda}.
\]
(Theorem~\ref{thm:asympt}, given that the characteristic polynomial $x-\lambda$ has root $\kappa=\lambda$)
\end{itemize}
In particular, the Conolly sequence~\cite{con} and the Tanny sequence~\cite{tanny} are asymptotic to $\frac{n}{2}$, as was known previously.

\subsection{Fibonacci Recurrence $a_n=a_{n-1}+a_{n-2}$}
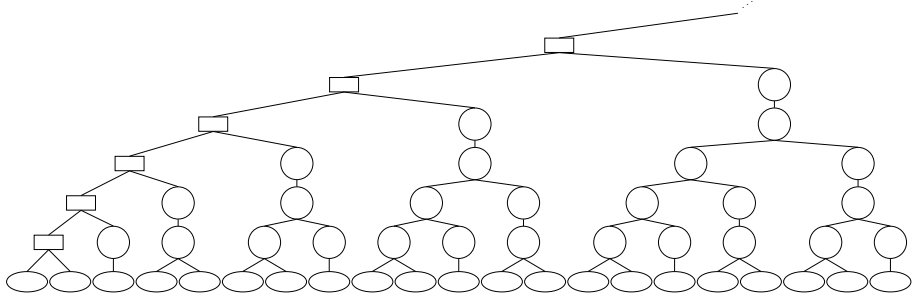
\begin{figure}
\begin{center}
\resizebox{\textwidth-0.5in}{!}{
%drawTree(levels = 6, a = [1, 1])
\tikzstyle{super}=[rectangle, draw]
\tikzstyle{regular}=[circle, draw]
\tikzstyle{leaf}=[ellipse, draw]
\begin{tikzpicture}[level distance = 1.5cm, grow=down]
\Tree[.{$\phantom{aaa}\iddots$}
[.\node[super]{$\phantom{0000}$};
  [.\node[super]{$\phantom{0000}$};
    [.\node[super]{$\phantom{0000}$};
      [.\node[super]{$\phantom{0000}$};
        [.\node[super]{$\phantom{0000}$};
          [.\node[super]{$\phantom{0000}$};
            [.\node[leaf]{$\phantom{0000}$};
            ]
            [.\node[leaf]{$\phantom{0000}$};
            ]
          ]
          [.\node[regular]{$\phantom{0000}$};
            [.\node[leaf]{$\phantom{0000}$};
            ]
          ]
        ]
        [.\node[regular]{$\phantom{0000}$};
          [.\node[regular]{$\phantom{0000}$};
            [.\node[leaf]{$\phantom{0000}$};
            ]
            [.\node[leaf]{$\phantom{0000}$};
            ]
          ]
        ]
      ]
      [.\node[regular]{$\phantom{0000}$};
        [.\node[regular]{$\phantom{0000}$};
          [.\node[regular]{$\phantom{0000}$};
            [.\node[leaf]{$\phantom{0000}$};
            ]
            [.\node[leaf]{$\phantom{0000}$};
            ]
          ]
          [.\node[regular]{$\phantom{0000}$};
            [.\node[leaf]{$\phantom{0000}$};
            ]
          ]
        ]
      ]
    ]
    [.\node[regular]{$\phantom{0000}$};
      [.\node[regular]{$\phantom{0000}$};
        [.\node[regular]{$\phantom{0000}$};
          [.\node[regular]{$\phantom{0000}$};
            [.\node[leaf]{$\phantom{0000}$};
            ]
            [.\node[leaf]{$\phantom{0000}$};
            ]
          ]
          [.\node[regular]{$\phantom{0000}$};
            [.\node[leaf]{$\phantom{0000}$};
            ]
          ]
        ]
        [.\node[regular]{$\phantom{0000}$};
          [.\node[regular]{$\phantom{0000}$};
            [.\node[leaf]{$\phantom{0000}$};
            ]
            [.\node[leaf]{$\phantom{0000}$};
            ]
          ]
        ]
      ]
    ]
  ]
  [.\node[regular]{$\phantom{0000}$};
    [.\node[regular]{$\phantom{0000}$};
      [.\node[regular]{$\phantom{0000}$};
        [.\node[regular]{$\phantom{0000}$};
          [.\node[regular]{$\phantom{0000}$};
            [.\node[leaf]{$\phantom{0000}$};
            ]
            [.\node[leaf]{$\phantom{0000}$};
            ]
          ]
          [.\node[regular]{$\phantom{0000}$};
            [.\node[leaf]{$\phantom{0000}$};
            ]
          ]
        ]
        [.\node[regular]{$\phantom{0000}$};
          [.\node[regular]{$\phantom{0000}$};
            [.\node[leaf]{$\phantom{0000}$};
            ]
            [.\node[leaf]{$\phantom{0000}$};
            ]
          ]
        ]
      ]
      [.\node[regular]{$\phantom{0000}$};
        [.\node[regular]{$\phantom{0000}$};
          [.\node[regular]{$\phantom{0000}$};
            [.\node[leaf]{$\phantom{0000}$};
            ]
            [.\node[leaf]{$\phantom{0000}$};
            ]
          ]
          [.\node[regular]{$\phantom{0000}$};
            [.\node[leaf]{$\phantom{0000}$};
            ]
          ]
        ]
      ]
    ]
  ]
]
\edge[draw=none]; {} ]
\end{tikzpicture}
}
\end{center}
\caption{Skeleton for the Fibonacci recurrence $\bk{1,1}$.}
\label{fig:fibskel}
\end{figure}

The first seven levels of the skeleton obtained from the Fibonacci recurrence are depicted in Figure~\ref{fig:fibskel}. The properties of the resulting trees/sequences/representations are described below.

\begin{itemize}
\item The tree $T_j$ is the family tree of $j$ generations of \quot{Fibonacci rabbits.} Classically, the Fibonacci sequence counts the number of pairs of rabbits present after $j$ generations where a mature pair produces a new pair every generation, and it takes one generation for a pair to mature. \typel{1} nodes correspond to immature pairs, and \typel{0} nodes correspond to mature pairs. (In this case, the process starts with a single, \emph{immature} pair.)
\item The skeleton $T$ is depicted in Figure~\ref{fig:fibskel}.
\item The $a$-Zeckendorf representation of $N$ is the Zeckendorf representation of $N$. (Proposition~\ref{prop:placeval})
\item The number of occurrences of the number $N$ in the sequence $L_{T,s}$ is one more than the number of zeroes at the end of the Zeckendorf representation of $N$, unless $N=F_{t+2}$ for some $t$, in which case the number of occurrences is $t+1+s$. When $s=0$, the resulting sequence is A316628 in the OEIS~\cite{oeis}. (Proposition~\ref{prop:ztreecount})
\item According to Theorem~\ref{thm:main}, the sequence $\seq{L_{T,s}\p{n}}_{n\geq1}$ satisfies the recurrence \[
L_{T,s}\p{n}=L_{T,s}\p{n-s-L_{T,s}\p{n-1}}+L_{T,s}\p{n-s-1-L_{T,s}\p{n-2}-L_{T,s}\p{n-s-2-L_{T,s}\p{n-2}}}.
\]
\item The sequence $\seq{L_{T,s}\p{n}}_{n\geq1}$ has the property that
\[
\lim_{n\to\infty}\frac{L_{T,s}\p{n}}{n}=1-\frac{1}{\phi}=0.381966\ldots,
\]
where $\phi=\frac{1+\sqrt{5}}{2}$.
(Theorem~\ref{thm:asympt}, given that the characteristic polynomial $x^2-x-1$ has positive root $\kappa=\phi$)
\end{itemize}

\subsection{Recurrence $a_n=a_{n-1}+2a_{n-2}+3a_{n-4}$}

\begin{figure}
\begin{center}
\begin{tabular}{c|c}
\resizebox{0.2in}{!}{
%drawTree(j = 0, levels = 0, a = [1, 2, 0, 3], height = 5)
\tikzstyle{super}=[rectangle, draw]
\tikzstyle{regular}=[circle, draw]
\tikzstyle{leaf}=[circle, draw]
\begin{tikzpicture}[level distance = 5.0cm, grow=down]
\Tree[.\node[leaf]{$\phantom{0000}$};
]

\end{tikzpicture}
} & 
\resizebox{0.2in}{!}{
%drawTree(j = 1, levels = 1, a = [1, 2, 0, 3], height = 5)
\tikzstyle{super}=[rectangle, draw]
\tikzstyle{regular}=[circle, draw]
\tikzstyle{leaf}=[circle, draw]
\begin{tikzpicture}[level distance = 5.0cm, grow=down]
\Tree[.\node[regular]{$\phantom{0000}$};
  [.\node[leaf]{$\phantom{0000}$};
  ]
]

\end{tikzpicture}
}\\\hline
\resizebox{0.5in}{!}{
%drawTree(j = 2, levels = 2, a = [1, 2, 0, 3], height = 5)
\tikzstyle{super}=[rectangle, draw]
\tikzstyle{regular}=[circle, draw]
\tikzstyle{leaf}=[circle, draw]
\begin{tikzpicture}[level distance = 5.0cm, grow=down]
\Tree[.\node[regular]{$\phantom{0000}$};
  [.\node[regular]{$\phantom{0000}$};
    [.\node[leaf]{$\phantom{0000}$};
    ]
    [.\node[leaf]{$\phantom{0000}$};
    ]
    [.\node[leaf]{$\phantom{0000}$};
    ]
  ]
]

\end{tikzpicture}
} & 
\resizebox{0.5in}{!}{
%drawTree(j = 3, levels = 3, a = [1, 2, 0, 3], height = 5)
\tikzstyle{super}=[rectangle, draw]
\tikzstyle{regular}=[circle, draw]
\tikzstyle{leaf}=[circle, draw]
\begin{tikzpicture}[level distance = 5.0cm, grow=down]
\Tree[.\node[regular]{$\phantom{0000}$};
  [.\node[regular]{$\phantom{0000}$};
    [.\node[regular]{$\phantom{0000}$};
      [.\node[leaf]{$\phantom{0000}$};
      ]
      [.\node[leaf]{$\phantom{0000}$};
      ]
      [.\node[leaf]{$\phantom{0000}$};
      ]
    ]
    [.\node[regular]{$\phantom{0000}$};
      [.\node[leaf]{$\phantom{0000}$};
      ]
    ]
    [.\node[regular]{$\phantom{0000}$};
      [.\node[leaf]{$\phantom{0000}$};
      ]
    ]
  ]
]

\end{tikzpicture}
}\\\hline
\resizebox{\textwidth*\real{0.35}-0.2in}{!}{
%drawTree(j = 4, a = [1, 2, 0, 3], height = 5)
\tikzstyle{super}=[rectangle, draw]
\tikzstyle{regular}=[circle, draw]
\tikzstyle{leaf}=[circle, draw]
\begin{tikzpicture}[level distance = 5.0cm, grow=down]
\Tree[.\node[regular]{$\phantom{0000}$};
  [.\node[regular]{$\phantom{0000}$};
    [.\node[regular]{$\phantom{0000}$};
      [.\node[regular]{$\phantom{0000}$};
        [.\node[leaf]{$\phantom{0000}$};
        ]
        [.\node[leaf]{$\phantom{0000}$};
        ]
        [.\node[leaf]{$\phantom{0000}$};
        ]
        [.\node[leaf]{$\phantom{0000}$};
        ]
        [.\node[leaf]{$\phantom{0000}$};
        ]
        [.\node[leaf]{$\phantom{0000}$};
        ]
      ]
      [.\node[regular]{$\phantom{0000}$};
        [.\node[leaf]{$\phantom{0000}$};
        ]
      ]
      [.\node[regular]{$\phantom{0000}$};
        [.\node[leaf]{$\phantom{0000}$};
        ]
      ]
    ]
    [.\node[regular]{$\phantom{0000}$};
      [.\node[regular]{$\phantom{0000}$};
        [.\node[leaf]{$\phantom{0000}$};
        ]
        [.\node[leaf]{$\phantom{0000}$};
        ]
        [.\node[leaf]{$\phantom{0000}$};
        ]
      ]
    ]
    [.\node[regular]{$\phantom{0000}$};
      [.\node[regular]{$\phantom{0000}$};
        [.\node[leaf]{$\phantom{0000}$};
        ]
        [.\node[leaf]{$\phantom{0000}$};
        ]
        [.\node[leaf]{$\phantom{0000}$};
        ]
      ]
    ]
  ]
]

\end{tikzpicture}
} & 
\resizebox{\textwidth*\real{0.5}-0.2in}{!}{
%drawTree(j = 5, levels = 5, a = [1, 2, 0, 3], height = 5)
\tikzstyle{super}=[rectangle, draw]
\tikzstyle{regular}=[circle, draw]
\tikzstyle{leaf}=[circle, draw]
\begin{tikzpicture}[level distance = 5.0cm, grow=down]
\Tree[.\node[regular]{$\phantom{0000}$};
  [.\node[regular]{$\phantom{0000}$};
    [.\node[regular]{$\phantom{0000}$};
      [.\node[regular]{$\phantom{0000}$};
        [.\node[regular]{$\phantom{0000}$};
          [.\node[leaf]{$\phantom{0000}$};
          ]
          [.\node[leaf]{$\phantom{0000}$};
          ]
          [.\node[leaf]{$\phantom{0000}$};
          ]
          [.\node[leaf]{$\phantom{0000}$};
          ]
          [.\node[leaf]{$\phantom{0000}$};
          ]
          [.\node[leaf]{$\phantom{0000}$};
          ]
        ]
        [.\node[regular]{$\phantom{0000}$};
          [.\node[leaf]{$\phantom{0000}$};
          ]
        ]
        [.\node[regular]{$\phantom{0000}$};
          [.\node[leaf]{$\phantom{0000}$};
          ]
        ]
        [.\node[regular]{$\phantom{0000}$};
          [.\node[leaf]{$\phantom{0000}$};
          ]
        ]
        [.\node[regular]{$\phantom{0000}$};
          [.\node[leaf]{$\phantom{0000}$};
          ]
        ]
        [.\node[regular]{$\phantom{0000}$};
          [.\node[leaf]{$\phantom{0000}$};
          ]
        ]
      ]
      [.\node[regular]{$\phantom{0000}$};
        [.\node[regular]{$\phantom{0000}$};
          [.\node[leaf]{$\phantom{0000}$};
          ]
          [.\node[leaf]{$\phantom{0000}$};
          ]
          [.\node[leaf]{$\phantom{0000}$};
          ]
        ]
      ]
      [.\node[regular]{$\phantom{0000}$};
        [.\node[regular]{$\phantom{0000}$};
          [.\node[leaf]{$\phantom{0000}$};
          ]
          [.\node[leaf]{$\phantom{0000}$};
          ]
          [.\node[leaf]{$\phantom{0000}$};
          ]
        ]
      ]
    ]
    [.\node[regular]{$\phantom{0000}$};
      [.\node[regular]{$\phantom{0000}$};
        [.\node[regular]{$\phantom{0000}$};
          [.\node[leaf]{$\phantom{0000}$};
          ]
          [.\node[leaf]{$\phantom{0000}$};
          ]
          [.\node[leaf]{$\phantom{0000}$};
          ]
        ]
        [.\node[regular]{$\phantom{0000}$};
          [.\node[leaf]{$\phantom{0000}$};
          ]
        ]
        [.\node[regular]{$\phantom{0000}$};
          [.\node[leaf]{$\phantom{0000}$};
          ]
        ]
      ]
    ]
    [.\node[regular]{$\phantom{0000}$};
      [.\node[regular]{$\phantom{0000}$};
        [.\node[regular]{$\phantom{0000}$};
          [.\node[leaf]{$\phantom{0000}$};
          ]
          [.\node[leaf]{$\phantom{0000}$};
          ]
          [.\node[leaf]{$\phantom{0000}$};
          ]
        ]
        [.\node[regular]{$\phantom{0000}$};
          [.\node[leaf]{$\phantom{0000}$};
          ]
        ]
        [.\node[regular]{$\phantom{0000}$};
          [.\node[leaf]{$\phantom{0000}$};
          ]
        ]
      ]
    ]
  ]
]

\end{tikzpicture}
}\\\hline
\multicolumn{2}{c}{\resizebox{\textwidth-0.2in}{!}{
%drawTree(j = 6, levels = 6, a = [1, 2, 0, 3], height = 5)
\tikzstyle{super}=[rectangle, draw]
\tikzstyle{regular}=[circle, draw]
\tikzstyle{leaf}=[circle, draw]
\begin{tikzpicture}[level distance = 5.0cm, grow=down]
\Tree[.\node[regular]{$\phantom{0000}$};
  [.\node[regular]{$\phantom{0000}$};
    [.\node[regular]{$\phantom{0000}$};
      [.\node[regular]{$\phantom{0000}$};
        [.\node[regular]{$\phantom{0000}$};
          [.\node[regular]{$\phantom{0000}$};
            [.\node[leaf]{$\phantom{0000}$};
            ]
            [.\node[leaf]{$\phantom{0000}$};
            ]
            [.\node[leaf]{$\phantom{0000}$};
            ]
            [.\node[leaf]{$\phantom{0000}$};
            ]
            [.\node[leaf]{$\phantom{0000}$};
            ]
            [.\node[leaf]{$\phantom{0000}$};
            ]
          ]
          [.\node[regular]{$\phantom{0000}$};
            [.\node[leaf]{$\phantom{0000}$};
            ]
          ]
          [.\node[regular]{$\phantom{0000}$};
            [.\node[leaf]{$\phantom{0000}$};
            ]
          ]
          [.\node[regular]{$\phantom{0000}$};
            [.\node[leaf]{$\phantom{0000}$};
            ]
          ]
          [.\node[regular]{$\phantom{0000}$};
            [.\node[leaf]{$\phantom{0000}$};
            ]
          ]
          [.\node[regular]{$\phantom{0000}$};
            [.\node[leaf]{$\phantom{0000}$};
            ]
          ]
        ]
        [.\node[regular]{$\phantom{0000}$};
          [.\node[regular]{$\phantom{0000}$};
            [.\node[leaf]{$\phantom{0000}$};
            ]
            [.\node[leaf]{$\phantom{0000}$};
            ]
            [.\node[leaf]{$\phantom{0000}$};
            ]
          ]
        ]
        [.\node[regular]{$\phantom{0000}$};
          [.\node[regular]{$\phantom{0000}$};
            [.\node[leaf]{$\phantom{0000}$};
            ]
            [.\node[leaf]{$\phantom{0000}$};
            ]
            [.\node[leaf]{$\phantom{0000}$};
            ]
          ]
        ]
        [.\node[regular]{$\phantom{0000}$};
          [.\node[regular]{$\phantom{0000}$};
            [.\node[leaf]{$\phantom{0000}$};
            ]
            [.\node[leaf]{$\phantom{0000}$};
            ]
            [.\node[leaf]{$\phantom{0000}$};
            ]
          ]
        ]
        [.\node[regular]{$\phantom{0000}$};
          [.\node[regular]{$\phantom{0000}$};
            [.\node[leaf]{$\phantom{0000}$};
            ]
            [.\node[leaf]{$\phantom{0000}$};
            ]
            [.\node[leaf]{$\phantom{0000}$};
            ]
          ]
        ]
        [.\node[regular]{$\phantom{0000}$};
          [.\node[regular]{$\phantom{0000}$};
            [.\node[leaf]{$\phantom{0000}$};
            ]
            [.\node[leaf]{$\phantom{0000}$};
            ]
            [.\node[leaf]{$\phantom{0000}$};
            ]
          ]
        ]
      ]
      [.\node[regular]{$\phantom{0000}$};
        [.\node[regular]{$\phantom{0000}$};
          [.\node[regular]{$\phantom{0000}$};
            [.\node[leaf]{$\phantom{0000}$};
            ]
            [.\node[leaf]{$\phantom{0000}$};
            ]
            [.\node[leaf]{$\phantom{0000}$};
            ]
          ]
          [.\node[regular]{$\phantom{0000}$};
            [.\node[leaf]{$\phantom{0000}$};
            ]
          ]
          [.\node[regular]{$\phantom{0000}$};
            [.\node[leaf]{$\phantom{0000}$};
            ]
          ]
        ]
      ]
      [.\node[regular]{$\phantom{0000}$};
        [.\node[regular]{$\phantom{0000}$};
          [.\node[regular]{$\phantom{0000}$};
            [.\node[leaf]{$\phantom{0000}$};
            ]
            [.\node[leaf]{$\phantom{0000}$};
            ]
            [.\node[leaf]{$\phantom{0000}$};
            ]
          ]
          [.\node[regular]{$\phantom{0000}$};
            [.\node[leaf]{$\phantom{0000}$};
            ]
          ]
          [.\node[regular]{$\phantom{0000}$};
            [.\node[leaf]{$\phantom{0000}$};
            ]
          ]
        ]
      ]
    ]
    [.\node[regular]{$\phantom{0000}$};
      [.\node[regular]{$\phantom{0000}$};
        [.\node[regular]{$\phantom{0000}$};
          [.\node[regular]{$\phantom{0000}$};
            [.\node[leaf]{$\phantom{0000}$};
            ]
            [.\node[leaf]{$\phantom{0000}$};
            ]
            [.\node[leaf]{$\phantom{0000}$};
            ]
            [.\node[leaf]{$\phantom{0000}$};
            ]
            [.\node[leaf]{$\phantom{0000}$};
            ]
            [.\node[leaf]{$\phantom{0000}$};
            ]
          ]
          [.\node[regular]{$\phantom{0000}$};
            [.\node[leaf]{$\phantom{0000}$};
            ]
          ]
          [.\node[regular]{$\phantom{0000}$};
            [.\node[leaf]{$\phantom{0000}$};
            ]
          ]
        ]
        [.\node[regular]{$\phantom{0000}$};
          [.\node[regular]{$\phantom{0000}$};
            [.\node[leaf]{$\phantom{0000}$};
            ]
            [.\node[leaf]{$\phantom{0000}$};
            ]
            [.\node[leaf]{$\phantom{0000}$};
            ]
          ]
        ]
        [.\node[regular]{$\phantom{0000}$};
          [.\node[regular]{$\phantom{0000}$};
            [.\node[leaf]{$\phantom{0000}$};
            ]
            [.\node[leaf]{$\phantom{0000}$};
            ]
            [.\node[leaf]{$\phantom{0000}$};
            ]
          ]
        ]
      ]
    ]
    [.\node[regular]{$\phantom{0000}$};
      [.\node[regular]{$\phantom{0000}$};
        [.\node[regular]{$\phantom{0000}$};
          [.\node[regular]{$\phantom{0000}$};
            [.\node[leaf]{$\phantom{0000}$};
            ]
            [.\node[leaf]{$\phantom{0000}$};
            ]
            [.\node[leaf]{$\phantom{0000}$};
            ]
            [.\node[leaf]{$\phantom{0000}$};
            ]
            [.\node[leaf]{$\phantom{0000}$};
            ]
            [.\node[leaf]{$\phantom{0000}$};
            ]
          ]
          [.\node[regular]{$\phantom{0000}$};
            [.\node[leaf]{$\phantom{0000}$};
            ]
          ]
          [.\node[regular]{$\phantom{0000}$};
            [.\node[leaf]{$\phantom{0000}$};
            ]
          ]
        ]
        [.\node[regular]{$\phantom{0000}$};
          [.\node[regular]{$\phantom{0000}$};
            [.\node[leaf]{$\phantom{0000}$};
            ]
            [.\node[leaf]{$\phantom{0000}$};
            ]
            [.\node[leaf]{$\phantom{0000}$};
            ]
          ]
        ]
        [.\node[regular]{$\phantom{0000}$};
          [.\node[regular]{$\phantom{0000}$};
            [.\node[leaf]{$\phantom{0000}$};
            ]
            [.\node[leaf]{$\phantom{0000}$};
            ]
            [.\node[leaf]{$\phantom{0000}$};
            ]
          ]
        ]
      ]
    ]
  ]
]

\end{tikzpicture}}
}\\
\end{tabular}
\end{center}
\caption{Trees $T_0$ through $T_6$ for the recurrence $\bk{1,2,0,3}$.}
\label{fig:atrees}
\end{figure}
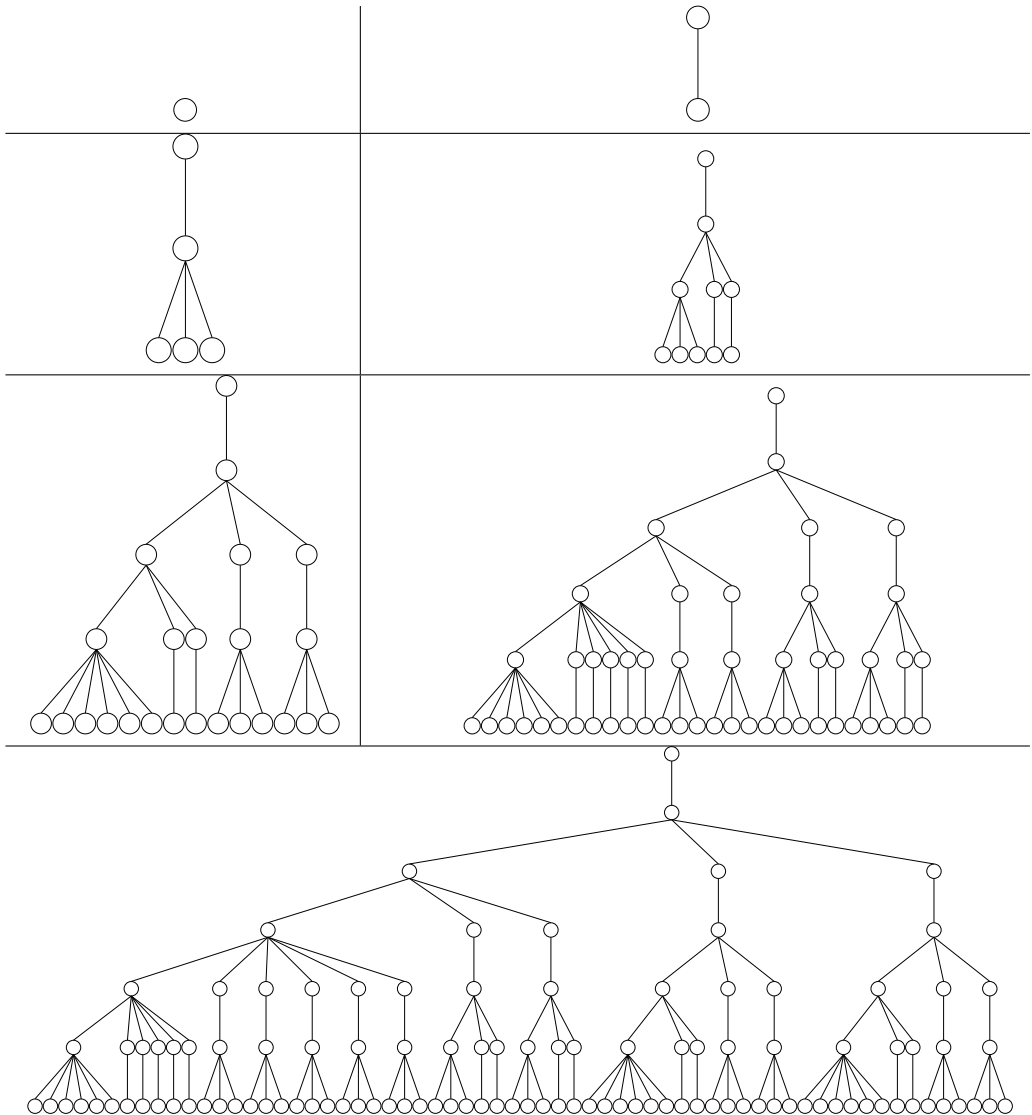

\begin{figure}
\begin{center}
\resizebox{\textwidth-0.5in}{!}{
\input{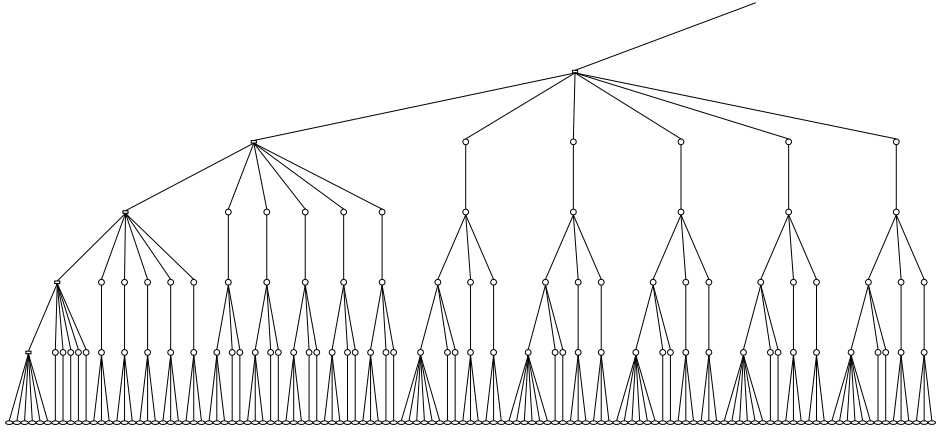}
}
\end{center}
\caption{Skeleton for the recurrence $\bk{1,2,0,3}$.}
\label{fig:askel}
\end{figure}

%The first seven levels of the skeleton obtained from the Fibonacci recurrence are depicted in Figure~\ref{fig:fibskel}. The properties of the resulting trees/sequences/representations are described below.
Here, we discuss the tree obtained from the recurrence $a=\bk{1,2,0,3}$ as an illustrative example that involves an order greater than $2$, a coefficient of $0$ in the recurrence, and an $a$-Zeckendorf representation system that is not a place value system. For this recurrence $\Lambda_1=1$, $\Lambda_2=3$, $\Lambda_3=3$, and $\Lambda_4=\Lambda=6$. Based on this, $\mu_0=1$, $\mu_1=2$, $\mu_2=2$, $\mu_3=4$, $\mu_4=4$, and $\mu_5=4$. The trees $T_0$ through $T_6$ are shown in Figure~\ref{fig:atrees}, and the skeleton is shown in Figure~\ref{fig:askel}. The first few terms of the sequence are $1,6,11,26,51,121,256,\ldots$.

\begin{itemize}
\item The $a$-valid digit sequences contain digits $0,1,2,3,4,5$. Every $1$ or $2$ is preceded by at least one $0$. Every $3$, $4$, or $5$ is preceded by at least three $0$'s. (Definition~\ref{def:valid})
\item The number of occurrences of the number $N$ in the sequence $\seq{L_{T,s}\p{n}}_{n\geq1}$ is one more than the number of zeroes at the end of the $a$-Zeckendorf representation of $N$, unless $N=a_t$ for some $t$, in which case the number of occurrences is $t+1+s$. (Proposition~\ref{prop:ztreecount})
\item According to Theorem~\ref{thm:main}, the sequence $\seq{L_{T,s}\p{n}}_{n\geq1}$ satisfies the recurrence 
\begin{align*}
L_{T,s}(n-L_{T,s}(n-1))&+L_{T,s}(n-1-L_{T,s}(n-2)-L_{T,s}(n-2-L_{T,s}(n-2)))\\
&+L_{T,s}(n-2-L_{T,s}(n-3)-L_{T,s}(n-3-L_{T,s}(n-3)))\\
&+L_{T,s}(n-3-L_{T,s}(n-4)-L_{T,s}(n-4-L_{T,s}(n-4))\\
&\phantom{+}-L_{T,s}(n-4-L_{T,s}(n-4)-L_{T,s}(n-4-L_{T,s}(n-4)))\\
&\phantom{+}-L_{T,s}(n-4-L_{T,s}(n-4)-L_{T,s}(n-4-L_{T,s}(n-4))\\
&\phantom{++}-L_{T,s}(n-4-L_{T,s}(n-4)-L_{T,s}(n-4-L_{T,s}(n-4)))))\\
&+L_{T,s}(n-4-L_{T,s}(n-5)-L_{T,s}(n-5-L_{T,s}(n-5))\\
&\phantom{+}-L_{T,s}(n-5-L_{T,s}(n-5)-L_{T,s}(n-5-L_{T,s}(n-5)))\\
&\phantom{+}-L_{T,s}(n-5-L_{T,s}(n-5)-L_{T,s}(n-5-L_{T,s}(n-5))\\
&\phantom{++}-L_{T,s}(n-5-L_{T,s}(n-5)-L_{T,s}(n-5-L_{T,s}(n-5)))))\\
&+L_{T,s}(n-5-L_{T,s}(n-6)-L_{T,s}(n-6-L_{T,s}(n-6))\\
&\phantom{+}-L_{T,s}(n-6-L_{T,s}(n-6)-L_{T,s}(n-6-L_{T,s}(n-6)))\\
&\phantom{+}-L_{T,s}(n-6-L_{T,s}(n-6)-L_{T,s}(n-6-L_{T,s}(n-6))\\
&\phantom{++}-L_{T,s}(n-6-L_{T,s}(n-6)-L_{T,s}(n-6-L_{T,s}(n-6))))).
\end{align*}
\item The characteristic polynomial $x^4-x^3-2x^2-3$ has a real root around $2.1949$. Hence, by Theorem~\ref{thm:asympt}, the sequence $L_{T,s}$ has the property that
\[
\lim_{n\to\infty}\frac{L_{T,s}\p{n}}{n}\approx1-\frac{1}{2.1949}\approx0.5444.
\]
\end{itemize}

\section{Algorithms}\label{s:algs}
This appendix is devoted to the study of efficient (linear in the length of the $a$-Zeckendorf representation; logarithmic in the number being represented) algorithms for converting between $a$-Zeckendorf representations and numbers.
%We are now ready to give our algorithms. We begin with the algorithm to go from $a$-Zeckendorf representations to numbers.  Given an $a$-valid word $w$ that is the $a$-Zeckendorf representation of some number $N$, by Proposition~\ref{prop:ztree} the representation describes a traversal down the skeleton $T$ for recurrence $a$ to a leaf. Then, $N$ is the number of non-leftmost leaves at or to the left of that leaf. Equivalently, $N$ is the number of leaves to the left of that leaf.

%Suppose that $w$ has length $m+1$. Since $w$ does not begin with $0$, the first step in the path is from the $m^{th}$ supernode to a regular node or leaf. To the left of this step is the left subtree of $T$ rooted at the $\p{m-1}^{st}$ supernode, which has $a_m$ leaves, by Corollary~\ref{cor:ztree}. 
Algorithm~\ref{alg:zn} works by 
calculating the formula in Proposition~\ref{prop:Nform}.
%starting from $0$ and accumulating the value of $N$ as $w$ is traversed from left to right. As we traverse $w$, we travel down $T$ until we reach a leaf at the end. There are two ways that leaves can appear to the left of the leaf reached along this path:
%On Line~\ref{ln:znquadam} of Algorithm~\ref{alg:znquad}, we initialize $N$ to $a_m$ to account for these leaves from the start. Aside from these $a_m$ leaves, there are two other ways leaves can appear to the left of the leaf reached along the path:
%\begin{itemize}
%\item On a step from a node $v$ to a non-\typel{0} child, all leaves in the subtree rooted at $v$'s \typel{0} child are to the left of the reached node.
%\item On a step from a node $v$ to a \typel{\ell} child where $\ell>0$, all leaves in the subtrees rooted at children of $v$ of type greater than $0$ and less than $\ell$ are to the left of the reached node.
%\end{itemize}
%Any step to a non-\typel{0} child is to the root of some $T_j$. Line~\ref{ln:znquadleft} accumulates the leaf counts of subtrees rooted at \typel{0} nodes, and Line~\ref{ln:znquadnotleft} accumulates the other leaf counts. Algorithm~\ref{alg:znquad} is written traversing the path from bottom up, though it could have also been implemented top-down.
%Aside from the initial accumulation of $a_m$, each subsequent accumulation of the first type involves adding a value $\NN\p{j,t}$ for some $0\leq t<j$. Accumulations of the second type involve adding values of the form $\NN\p{j,j}$.

\begin{algorithm}
\caption{Algorithm for converting an $a$-Zeckendorf representation to a number}
\label{alg:zn}
\begin{algorithmic}[1]
\Procedure {ZeckToNum}{$a$; $\dseq{d}{M}$}\Comment{Inputs: a recurrence $a$ and an $a$-valid digit sequence}
%\State $m\gets\text{length}(w)-1$
\State $N\gets 0$\Comment{The value to be accumulated and eventually returned}\label{ln:znam}
%\State $lastNonzeroIndex\gets undefined$\Comment{Will track index of most recent nonzero digit}
\State $j\gets-1$\Comment{Most recent nonzero index; $-1$ placeholder initially}
\For{$t$ from $M$ down to $0$}\Comment{Loop through digits from left to right}
\If {$d_t\neq0$}
\If {$j=-1$}
\State $N\gets N + a_M$\Comment{Corresponds to initial left subtree}
\Else
\State $N\gets N + \NN\p{j,t}$\Comment{Corresponds to internal left subtree}\label{ln:znleft}
\EndIf
\State $j\gets t$\Comment{Track index of most recent nonzero digit}
\State $N\gets N + \pb{d_t-1}\NN\p{t,t}$\Comment{Corresponds to non-left subtrees}\label{ln:znnotleft}
\EndIf
\EndFor
\State \textbf{return} $N$\Comment{The value is $N$}
\EndProcedure
\end{algorithmic}
\end{algorithm}

Algorithm~\ref{alg:zn} as stated requires computation of the values $\NN\p{j,t}$. Using the formula from Proposition~\ref{prop:tjleafcount} for $\NN\p{j,t}$, we can compute $\NN\p{j,t}$ in $O\p{j-t}=O\p{M}$ time. These values appear inside a loop, which makes this version of the algorithm quadratic. We can reduce the time complexity to linear by exploiting the recurrence in Lemma~\ref{lem:Nrec} to update the values $\NN\p{j,t}$ as we go, thereby reducing the overall runtime to $O\p{M}=O\p{\log N}$.  This is done in Algorithm~\ref{alg:znfast}.

\begin{algorithm}
\caption{Linear algorithm for converting an $a$-Zeckendorf representation to a number}
\label{alg:znfast}
\begin{algorithmic}[1]
\Procedure {ZeckToNumFast}{$a$; $\dseq{d}{M}$}\Comment{Inputs: a recurrence $a$ and an $a$-valid digit sequence}
%\State $m\gets\text{length}(w)-1$
\State $N\gets 0$\Comment{The value to be accumulated and eventually returned}\label{ln:znfastam}
%\State $lastNonzeroIndex\gets undefined$\Comment{Will track index of most recent nonzero digit}
\State $left\gets a_M$\Comment{Current leaf count of left subtree}
\State $j\gets-1$\Comment{Most recent nonzero index; $-1$ placeholder initially}
\For{$t$ from $M$ down to $0$}\Comment{Loop through digits from left to right}
\If {$d_t\neq0$}
\State $N\gets N + left$\Comment{Corresponds to left subtree}\label{ln:znfastleft}
\State $left\gets \frac{1}{\Lambda-1}\pb{a_{t+1}-\Lambda_1a_t+\pb{\Lambda_1-1}a_{t-1}}$\Comment{Proposition~\ref{prop:tjleafcount} formula for $\NN\p{t,t-1}$, the leaf count of the next left subtree}
\State $j\gets t$\Comment{Track index of most recent nonzero digit}
\State $N\gets N + \pb{d_t-1}\NN\p{t,t}$\Comment{Corresponds to non-left subtrees}\label{ln:znfastnotleft}
\Else\Comment{Zero digit; still need to update $left$}
\State $left\gets left-\pb{\Lambda_{j-t+1}-1}\frac{a_t-a_{t-1}}{\Lambda-1}$\Comment{Lemma~\ref{lem:Nrec} to reduce $\NN\p{j,t}$ to $\NN\p{j,t-1}$.}
\EndIf
\EndFor
\State \textbf{return} $N$\Comment{The value is $N$}
\EndProcedure
\end{algorithmic}
\end{algorithm}

Our algorithm for taking a number $N$ and producing its $a$-Zeckendorf representation is quite similar. It is also reminiscent of the greedy algorithm for calculating Zeckendorf representations. Suppose $N$ is such that $a_M\leq N<a_{M+1}$. Then, the $a$-Zeckendorf representation of $N$ has $M+1$ digits. Like the previous algorithm, we work left to right through the digits, accumulating a value that starts from $0$. In this case, we accumulate the value of the number being represented so far. For each index, we first determine which digits could validly appear there. If there is a nonzero digit we could put there while making the value accumulated not exceed $N$, we put the largest such value there and increase the accumulated value accordingly. Otherwise, we make that digit $0$. The process is implemented in Algorithm~\ref{alg:nz}

\begin{algorithm}
\caption{Algorithm for converting a number to an $a$-Zeckendorf representation}
\label{alg:nz}
\begin{algorithmic}[1]
\Procedure {NumToZeck}{$a$, $N$}\Comment{Inputs: $a$ a recurrence; $N$ a positive integer}
\State $M$ is the unique value such that $a_M\leq N<a_{M+1}$
\State $sofar\gets 0$\Comment{The value of our representation so far}
%\State $lastNonzeroIndex\gets undefined$\Comment{Will track index of most recent nonzero digit}
\State $j\gets-1$\Comment{Most recent nonzero index; $-1$ placeholder initially}
\For{$t$ from $M$ down to $0$}\Comment{Loop through digits from left to right}
\If {$j=-1$ or $N-sofar\geq\NN\p{j,t}$}\Comment{We have a nonzero digit}\label{ln:nznonzero}
\State $d\gets\fl{\frac{n-sofar-\NN\p{j,t}}{\NN\p{t,t}}}$\Comment{How many copies of $T_t$ are to the left?}\label{ln:nzcopies}
\State $sofar\gets sofar+d\cdot\NN\p{t,t}$\Comment{Accumulate the leaves from the $T_t$ copies}
\State $d_t\gets d+1$\Comment{Digit is one more than number of $T_t$ copies}
\State $sofar\gets sofar+\NN\p{j,t}$\Comment{Accumulate the leaves from the left subtree}
\State $j\gets t$\Comment{Track index of most recent nonzero digit}
\Else
\State $d_t\gets0$\Comment{We have a zero digit}
\EndIf
\EndFor
\State \textbf{return} $\dseq{d}{M}$\Comment{The representation of $N$}
\EndProcedure
\end{algorithmic}
\end{algorithm}

To check if we have a nonzero digit, we actually check whether the difference between $N$ and the value represented so far is at least as large as the number of leaves in the left subtree. If it is, all of those nodes must be to the left of the leaf we eventually reach, meaning we have a nonzero digit. We then need to determine which copy of $T_t$ below the current node the path goes to. Each copy has $\NN\p{t,t}$ leaves, so the computation on Line~\ref{ln:nzcopies} accomplishes this task. The rest of the loop is bookkeeping.

Like Algorithm~\ref{alg:zn}, Algorithm~\ref{alg:nz}, as given, has quadratic runtime in $M$. Again, we can use the recurrence from Lemma~\ref{lem:Nrec} to reduce this runtime to linear. The result is Algorithm~\ref{alg:nzfast}.

\begin{algorithm}
\caption{Linear algorithm for converting a number to an $a$-Zeckendorf representation}
\label{alg:nzfast}
\begin{algorithmic}[1]
\Procedure {NumToZeckFast}{$a$, $N$}\Comment{Inputs: $a$ a recurrence; $N$ a positive integer}
\State $M$ is the unique value such that $a_M\leq N<a_{M+1}$
\State $sofar\gets 0$\Comment{The value of our representation so far}
%\State $lastNonzeroIndex\gets undefined$\Comment{Will track index of most recent nonzero digit}
\State $left\gets a_M$\Comment{Current leaf count of left subtree}
\State $j\gets-1$\Comment{Most recent nonzero index; $-1$ placeholder initially}
\For{$t$ from $M$ down to $0$}\Comment{Loop through digits from left to right}
\If {$N-sofar\geq left$}\Comment{We have a nonzero digit}
\State $Ltt\gets\frac{a_{t+1}-a_t}{\Lambda-1}$\Comment{Store value of $\NN\p{t,t}$ for readability}
\State $d\gets\fl{\frac{n-sofar-left}{Ltt}}$\Comment{How many copies of $T_t$ are to the left?}
\State $sofar\gets sofar+d\cdot Ltt$\Comment{Accumulate the leaves from the $T_t$ copies}
\State $d_t\gets d+1$\Comment{Digit is one more than number of $T_t$ copies}
\State $sofar\gets sofar+left$\Comment{Accumulate the leaves from the left subtree}
\State $left\gets \frac{1}{\Lambda-1}\pb{a_{t+1}-\Lambda_1a_t+\pb{\Lambda_1-1}a_{t-1}}$\Comment{Proposition~\ref{prop:tjleafcount} formula for $\NN\p{t,t-1}$, the leaf count of the next left subtree}
\State $j\gets t$\Comment{Track index of most recent nonzero digit}
\Else\Comment{Zero digit; still need to update $left$}
\State $d_t\gets0$\Comment{We have a zero digit}
\State $left\gets left-\pb{\Lambda_{j-t+1}-1}\frac{a_t-a_{t-1}}{\Lambda-1}$\Comment{Lemma~\ref{lem:Nrec} to reduce $\NN\p{j,t}$ to $\NN\p{j,t-1}$.}
\EndIf
\EndFor
\State \textbf{return} $\dseq{d}{M}$\Comment{The representation of $N$}
\EndProcedure
\end{algorithmic}
\end{algorithm}

\end{document}